\documentclass[12pt,a4paper]{amsart}

\usepackage{amssymb,amscd,amsmath,mathrsfs}

\input xy
\xyoption{all}

\newif\ifEU\IfFileExists{euscript.sty}{\EUtrue\usepackage[mathcal]{euscript}}{}

\makeatletter
\renewcommand{\phi}{\varphi}
\renewcommand{\to}{\longrightarrow}
\renewcommand{\le}{\leqslant}
\renewcommand{\geq}{\geqslant}
\renewcommand{\epsilon}{\varepsilon}

\renewcommand{\kappa}{\varkappa}
\DeclareMathOperator{\CAT}{\mathbf{CAT}}
\DeclareMathOperator{\Holim}{Holim}

\DeclareMathOperator{\cone}{cone}
\DeclareMathOperator{\diag}{diag}

 \DeclareMathOperator{\Hom}{Hom}
 
\DeclareMathOperator{\id}{id} \DeclareMathOperator{\Cat}{{\bf
Cat}} \DeclareMathOperator{\Mor}{Mor}
 \DeclareMathOperator{\Ho}{Ho}

 \DeclareMathOperator{\im}{Im}

 \DeclareMathOperator{\Ar}{Ar}

\DeclareMathOperator{\Ob}{Ob}
\newcommand{\la}{\textrm{\huge${}_{\ulcorner}$}}
\newcommand{\ra}{\textrm{\LARGE${\lrcorner}$}}

\newcommand {\holp}{\underrightarrow{\Holim}}
\newcommand {\holo}{\underleftarrow{\Holim}}

\newcommand{\lra}[1]{\bl{#1}\longrightarrow\relax}
\newcommand{\bl}[1]{\buildrel #1\over}

\newcommand{\cc}{\mathscr}
\newcommand{\bb}{\mathbf}
\newcommand{\ps}{\oplus}

\newcommand{\iso}{\simeq}
\newcommand{\op}{{\textrm{\rm op}}}

\newcommand{\wt}{\widetilde}

\newcommand{\es}[3]{0\to #1\to #2\to #3\to 0}
\newcommand{\ifff}{if and only if }

\newcommand {\dia}{{\mathcal Dia^{\star}}}
\newcommand {\di}{\mathcal Dia}
\newcommand {\ord}{\mathcal Ord}
\newcommand {\ordd}{\mathcal Ord^{\star}}
\newcommand {\dirf}{\mathcal{D}\mathit{irf}}
\newcommand{\homa}[2]{\Hom(#1,#2)}

\newtheorem{thm}{\sc Theorem}[section]
\newtheorem{prop}[thm]{\sc Proposition}
\newtheorem*{funcax}{\bf\em Functoriality Axiom}
\newtheorem*{isomax}{\bf\em Isomorphism Axiom}
\newtheorem*{disjax}{\bf\em Disjoint Union Axiom}
\newtheorem*{kanax}{\bf\em Homotopy Kan Extension Axioms}
\newtheorem*{baseax}{\bf\em Base Change Axiom}
\newtheorem{cor}[thm]{\sc Corollary}
\newtheorem{lem}[thm]{\sc Lemma}

\newtheorem*{question}{\sc Question}

\newtheorem*{observation}{\sc Observation}
\newtheorem*{immersion}{\sc The Immersion Lemma}
\newtheorem*{rem}{\sc Remark}
\newtheorem*{exs}{\sc Examples}

\newtheorem*{ex}{\sc Example}

\newtheorem*{theo}{\sc Theorem}

\newtheorem*{defs}{\sc Definition}
\newtheorem*{conv}{\sc Convention}

\makeatother

\begin{document}

\footskip30pt


\title{Systems of Diagram Categories and K-theory. I}
\author{Grigory Garkusha}
\address{International Centre for Theoretical Physics,
Strada Costiera 11, I-34014, Trieste, Italy}
\urladdr{www.ictp.trieste.it/$\sim$garkusha}
\email{ggarkusha@mail.ru}
\thanks{Supported by the ICTP Research Fellowship\\ Russian version is available at author's homepage}
\begin{abstract}
To any left system of diagram categories or to any left pointed
d\'e\-ri\-va\-teur a $K$-theory space is associated. This
$K$-theory space is shown to be canonically an infinite loop space
and to have a lot of common properties with Waldhausen's
$K$-theory. A weaker version of additivity is shown. Also,
Quillen's $K$-theory of a large class of exact categories
including the abelian categories is proved to be a retract of the
$K$-theory of the associated d\'erivateur.
\end{abstract}
\maketitle

\thispagestyle{empty}
\pagestyle{plain}

\tableofcontents

\section*{Introduction}

The object of this paper is to construct the Waldhausen $K$-theory
for systems of diagram categories and pointed d\'erivateurs
respectively. The general formalism for them was developed by
Heller~\cite{H}, Grothendieck~\cite{G}, and Franke~\cite{F}.

A system of diagram categories (respectively d\'erivateur) is a
hyperfunctor
   $$\bb B:\di^{\op}\lra{}\CAT$$
defined on a 2-subcategory of the 2-category of small categories
$\Cat$ (we shall refer to $\di$ as a {\it category of diagrams})
and taking values in the ``large" 2-category of categories $\CAT$.
We also require $\di$ to contain the 2-category of posets (=finite
partially ordered sets) $\ord$. A typical example of such a
hyperfunctor is given by the map
   $$I\in\di\longmapsto\Ho\cc C^I,$$
where $\cc C$ is a closed model category, $\Ho\cc C^I$ is the
homotopy category of the functor category $\cc C^I$ (the structure
of a closed model category on $\cc C^I$ is naturally defined for
appropriate diagrams).

First we define the $S.$-construction for an appropriate $\bb B$
and then its $K$-theory space $K(\bb B)$ as the loop space of
$|i.S.\bb B|$, where $iS_n\bb B$ stands for the subcategory of
isomorphisms in each category $S_n\bb B$, $n\geq 0$. The space
$K(\bb B)$ is canonically an infinite loop space by Segal's
machine~\cite{S}. The additivity theorem is discussed for this
$K$-theory. Developing Waldhausen's machinery~\cite{W} for objects
in question the additivity theorem implies that one can also think
of $K(\bb B)$ in terms of the following connective
$\Omega$-spectrum. Namely, it is given by the sequence of spaces
   $$\Omega|i.S.\bb B|,\Omega|i.S.S.\bb B|,\ldots,\Omega|i.S.^n\bb B|,\ldots$$
where the multisimplicial objects $i.S.^n\bb B$, $n\geq 1$, are
obtained by iterating the $S.$-con\-struc\-tion. Though the
additivity theorem remains open in the general case (see
also~\cite[Conjecture~3]{M}), a weaker version does hold.

\begin{theo}
The additivity theorem is valid for the space
   $$\Omega^\infty|i.S.^\infty\bb B|= \lim_n\Omega^n|i.S.^n\bb B|$$
excluding pathological cases we never have in practice.
\end{theo}

The strong form of additivity is shown in~\cite{Gar} for
complicial d\'erivateurs.

We also remark that the Grothendieck group $K_0(\cc E)$ of an
exact category $\cc E$ is naturally isomorphic to the group
$K_0(\bb D^b(\cc E))$ of the associated d\'erivateur $\bb D^b(\cc
E)$ given by the hyperfunctor
   $$I\longmapsto D^b(\cc E^I)$$
where $D^b(\cc E^I)$ is the derived category for the exact functor category $\cc E^I$.

One can obtain some relation between Quillen's $K$-theory $K(\cc
E)$ and $K(\bb D^b(\cc E))$ for a large class of exact categories
including the abelian categories.

\begin{theo}
Let $\cc E$ be an extension closed full exact subcategory of an
abelian category $\cc A$ satisfying the conditions of the
Resolution Theorem. That is

$(1)$ if $\es {M'}M{M''}$ is exact in $\cc A$ and $M,M''\in\cc E$,
then $M'\in\cc E$ and

$(2)$ for any object $M\in\cc A$ there is a finite resolution
$0\to P_n\to P_{n-1}\to\cdots\to P_0\to M\to 0$ with $P_i\in\cc
E$.

Then a natural map
   $$K(\rho):K(\cc E)\to K(\bb D^b(\cc E))$$
is a split inclusion in homotopy. There is a map
   $$p:K(\bb D^b(\cc E))\to K(\cc E)$$
which is left inverse to it. That is $p\circ K(\rho)$ is homotopic
to the identity. In particular, each $K$-group $K_n(\cc E)$ is a
direct summand of $K_n(\bb D^b(\cc E))$.
\end{theo}

The problem whether Quillen's $K$-theory $K(\cc E)$ can be
reconstructed from $\bb D^b(\cc E)$ (``the first Maltsiniotis
conjecture'') remains open.

I would like to thank Denis-Charles Cisinski and Amnon Neeman for
helpful e-discussions.

\section{Systems of diagram categories}\label{fgh}

In this section the reader will be attacked by categorical
formalities and definitions of concepts. A lot of analogous
statements for this section are also stated in~\cite{F}. Since the
latter work contains a lot of slight errata and this paper is one
of the first in this direction, we prove them once more in details
{\it in order to be sure\/} that nothing goes wrong. This makes
therefore our paper self-contained. We follow here the original
terminology of Franke~\cite{F}. The author must admit that the
formal portion of this section is maximally shortened.

\subsection{Notations} Let $I$ be a category. For a subcategory $J$ of $I$
and $x\in I$, we shall denote by $J/x$ the following comma category.
Objects are the pairs $(y,\phi)$ where $y\in J$ and $\phi:y\to x$ is a morphism
in $I$. Morphisms from $(y,\phi)$ to $(y',\phi')$ are given by morphisms
$\psi:y\to y'$ in $J$ such that $\phi=\phi'\psi$. The category $J\setminus x$
consists of pairs $(y,\phi)$ with $y\in J$ and $\phi:x\to y$. Morphisms
are defined similar to those of $J/x$. If $K\subseteq\Ob I$ is a subclass
of objects, we shall denote by $I-K$ the full subcategory of $I$ with
the class of objects $I-K$. In particular, if $K=\{x\}$ has just one
object, we shall also denote this subcategory by $I-x$. If $f:J\to I$
is a functor, the categories $f/x$ and $f\setminus x$ have objects
$(y\in J,\phi:f(y)\to x)$ and $(y\in J,\phi:x\to f(y))$. If $f$ is the inclusion
of a subcategory, this is the same as $J/x$ and $J\setminus x$.

Given a non-negative integer $n$, by $\Delta^n$ denote the totally
ordered set $\{0<1<\cdots<n\}$. For $i\le n+1$, the map
$d_i:\Delta^n\to\Delta^{n+1}$ is the monotonic injection not
containing $i$ in its image and $s_i:\Delta^n\to\Delta^{n-1}$ is
the monotonic surjection satisfying $s_i(i)=s_i(i+1)$.

\subsection{The axioms}

For the notions of the 2-category and 2-functor we refer the reader to~\cite{Mac}.
In what follows we use the term ``poset" as an abbreviation of ``finite
partially ordered set". Every poset can be considered as a category in
which $\homa xy$ has precisely one element $x\le y$, and is empty otherwise.
The 2-category of the posets (respectively the finite categories without cycles)
we shall denote by $\ord$ (respectively by $\dirf$).

Let $\di$ be a full 2-subcategory of the 2-category $\Cat$ of
small categories that contains the 2-category $\ord$. In what
follows we assume that $\di$ satisfies the following conditions:

\begin{enumerate}
\item $\di$ is closed under finite sums and finite products;
\item for any functor $f:I\to J$ in $\di$ and for any object $y$ of $J$,
      the categories $f/y$ and $f\setminus y$ are in $\di$.
\end{enumerate}

We shall also refer to $\di$ as a {\it category of diagrams}.

Given $I\in\di$, let $I^\star$ be $I$ with an initial and final object $\star$ added.
For any $x$ and $y$ in $I$, there is the unique morphism from $x$ to $y$
in $I^\star$ which factorizes through $\star$. We shall refer to this
morphism as the zero morphism. If $I\in\ord$ and $x\le y$ there is one more morphism from $x$
to $y$ in $I^\star$, and there are no other morphisms. The composition is defined
in the obvious way. Let $\dia$ be a 2-subcategory of the 2-category $\Cat$ whose
objects are those of $\di$ and
whose horizontal morphisms $I\to J$ are given by functors
$I^\star\to J^\star$ mapping $\star$ to $\star$, and let bimorphisms be
natural transformations between functors from $I^\star$ to $J^\star$.

A {\it presystem of diagram categories of the domain $\di$\/} or just a
{\it presystem of diagram categories\/} is a functor
   \begin{equation}\label{dia}
    \bb C:\dia^{\op}\to\CAT
   \end{equation}
from $\dia$ to the category $\CAT$ of categories (not necessarily small) satisfying
the Functoriality Axiom below. So to each category $I$ in $\dia$ there is associated
a category $\bb C_I$, and to each map $f:I\to J$ in $\dia$ a functor
$f^*=\bb C(f):\bb C_J\to\bb C_I$.

\begin{funcax}{\rm The following conditions hold:

\begin{itemize}

\item[$\diamond$] to each natural transformation $\phi:f\to g$ a
natural transformation $\phi^*:f^*\to g^*$ is associated and the
maps $f\to f^*$ and $\phi\to\phi^*$ define a functor from $\homa
IJ$ to the category of functors from $\bb C_J$ to $\bb C_I$;

\item[$\diamond$] if
   $$\xymatrix{K\ar[r]^f &I\ar@/^/[r]^g \ar@/_/[r]_{g'} &J\ar[r]^h &L}$$
are morphisms and $\phi:g\to g'$ is a bimorphism, then
$f^*\circ\phi^*=(\phi\circ f)^*$ and $\phi^*\circ h^*=(h\circ\phi)^*$.

\end{itemize}
}\end{funcax}

From now on let us fix a category of diagrams $\di$.
To any category $\cc C$ one associates a presystem of diagram categories
which takes a category $I$ of $\dia$ to the functor category
   $$\cc C^{I^\star}=\homa {I^\star}{\cc C}$$
and a map $f:I\to J$ to the map
   $$f^*:\cc C^{J^\star}\to\cc C^{I^\star},\ \ \ X\longmapsto X\circ f.$$

A {\it morphism\/} $F:\bb C\to\bb C'$ between two presystems of diagram categories
$\bb C$ and $\bb C'$ consists of the following data:
\begin{enumerate}
\item for any $I\in\dia$, a functor $F:\bb C_I\to\bb C'_I$;
\item for any map $f:I\to J$ in $\dia$, an isomorphism of functors
$\iota_{F,f}:f^*F\lra{\sim}Ff^*$.
\end{enumerate}

We also assume the following conditions to hold for $\iota_{F,f}$:

\begin{itemize}
\item[$\diamond$] for any $I\in\dia$, $\iota_{F,1_I}=1_F$;
\item[$\diamond$] for any two composable maps $I\lra{f}J\lra{g}K$ in $\dia$, the diagram
   $$\xymatrix{f^*g^*F\ar[rr]^{\iota_{F,gf}} \ar[dr]_{f^*\iota_{F,g}} && Ff^*g^* & \\
               &f^*Fg^*\ar[ur]_{\iota_{F,f}g^*}&}$$
is commutative; \item[$\diamond$] for any bimorphism $\phi:f\to g$
in $\dia$, we have the following commutative square.
   $$\xymatrix{f^*F\ar[r]^{\iota_{F,f}}\ar[d]_{\phi^*F} &Ff^*\ar[d]^{F\phi^*}\\
               g^*F\ar[r]_{\iota_{F,g}}&Fg^*}$$
\end{itemize}
A morphism $F:\bb C\to\bb C'$ is an {\it equivalence\/} if for any $I\in\dia$ the
functor $F:\bb C_I\to\bb C'_I$ is an equivalence of categories.

Let $\bb A,\bb B,\bb C$ be presystems of diagram categories. By the {\it fibred product\/}
of a pair of morphisms $F:\bb A\to\bb C$ and $G:\bb B\to\bb C$ is meant the following data:

\begin{itemize}
\item[$\diamond$] for any $I\in\dia$, the category $\prod(F,G)_I$ whose objects are the
triples
   $$(A,c,B),\ \ \ A\in\bb A_I,\ \ \ B\in\bb B_I,\ \ \ c:F(A)\lra{\sim}G(B),$$
and where a morphism from $(A,c,B)$ to $(A',c',B')$ is a pair of morphisms
$(a,b)$ compatible with the isomorphisms $c$ and $c'$;

\item[$\diamond$] for any map $f:I\to J$ in $\dia$, the functor
   $$f^*=f^*_{\prod(F,G)}:\prod(F,G)_J\to\prod(F,G)_I$$
defined by
   $$(A,c,B)\longmapsto(f^*_{\bb A}(A),\iota_{G,f}\circ f^*_{\bb C}(c)\circ\iota^{-1}_{F,f},f^*_{\bb B}(B)).$$
\end{itemize}

\begin{prop}\label{qqq}
The above data determine a presystem of diagram categories
   \begin{equation}\label{dff}
    \prod(F,G):\dia^{op}\to\CAT.
   \end{equation}
\end{prop}

\begin{proof}
Let us show that~\eqref{dff} is a functor.
For this, consider two composable maps $I\lra{g}J\lra{f}K$ in $\dia$. We have
   \begin{align*}
   \iota_{G,g}f^*\circ g^*(\iota_{G,f}f^*(c)\iota^{-1}_{F,f})\circ\iota^{-1}_{F,g}f^*&=\\
     \underbrace{
     \iota_{G,g}f^*\circ g^*(\iota_{G,f})}_{\iota_{G,fg}}
     \circ g^*(f^*(c))\circ&\overbrace{g^*(\iota^{-1}_{F,f})\circ\iota^{-1}_{F,g}f^*}^{\iota^{-1}_{F,fg}}
     =\iota_{G,fg}(fg)^*(c)\iota^{-1}_{F,fg}.
   \end{align*}
We see that $(fg)^*=g^*f^*:\prod(F,G)_K\to\prod(F,G)_I$.

Since the diagram
   $$\xymatrix{
     Ff^*(A)\ar[r]^{\iota^{-1}_{F,f}}\ar[d]_{F\phi^*} & f^*F(A)\ar[r]^{f^*(c)}\ar[d]_{\phi^*F}
     &f^*G(B)\ar[r]^{\iota_{G,f}}\ar[d]^{\phi^*G} &Gf^*(B)\ar[d]^{G\phi^*}\\
     Fg^*(A)\ar[r]_{\iota^{-1}_{F,g}} & g^*F(A)\ar[r]_{g^*(c)} & g^*G(B)\ar[r]_{\iota_{G,g}}&Gg^*(B)}$$
is commutative for any morphisms $f,g:I\to J$ and any bimorphism $\phi:f\to g$ in $\dia$, the map
   $$(f^*(A),\iota_{G,f}f^*(c)\iota^{-1}_{F,f},f^*(B))\longmapsto
     (g^*(A),\iota_{G,g}g^*(c)\iota^{-1}_{F,g},g^*(B))$$
yields a map $\phi^*$ between $f^*,g^*:\prod(F,G)_J\to\prod(F,G)_I$. The functoriality axiom
is directly verified and left to the reader.
\end{proof}

Let $\Delta^n=\{0<\cdots<n\}\in\di$. If there is no liklihood
of confusion, we also denote $\Delta^0$ by 0. Given $I\in\di$ and $x\in I$, let $i_{x,I}:0\to I$ be the
functor sending 0 to $x$. For $A\in\bb C_I$ let $A_x=i_{x,I}^*A$.

Given $I\in\di$ there is a natural functor
   $$dia_I:\bb C_I\to\Hom(I,\bb C_0).$$
It is constructed as follows. For any $x\in I$ we put
$dia_I(B)(x)=B_x$. Every morphism $\alpha:x\to y$ in $I$ yields a
natural transformation $\alpha:i_{x,I}\to i_{y,I}$. Then
$dia_I(B)(\alpha):=\alpha^*:B_x\to B_y$.

Let us consider the following axioms listed below.

\begin{isomax}{\rm A morphism $f:A\to B$ in $\bb C_I$
is an isomorphism iff $dia_I(f)$ is so in $\Hom(I,\bb C_0)$. In
other words, it is an isomorphism iff $f_x:A_x\to B_x$ is so for
all $x\in I$. }\end{isomax}

\begin{disjax}{\rm
(a) If $I=I_1\coprod I_2$ is a disjoint union of its full
subcategories $I_1$ and $I_2$, then the inclusions
$i_{1;2}:I_{1;2}\to I$ define an equivalence of categories
   $$(i^*_1,i^*_2):\bb C_I\lra{\sim}\bb C_{I_1}\times\bb C_{I_2}.$$

(b) $\bb C_\emptyset$ is a trivial category (having precisely one
morphism between any pair of objects).
}\end{disjax}

\begin{kanax}{\rm
The left homotopy Kan extension axiom requires that for any
functor $f:I\to J$, the functor $f^*:\bb C_J\to\bb C_I$ has a left
adjoint $f_!:\bb C_I\to\bb C_J$. By symmetry, the right homotopy
Kan extension axiom says that $f^*$ has a right adjoint $f_*:\bb
C_I\to\bb C_J$. Below we shall also refer to the functors $f_!$
and $f_*$ as left and right homotopy Kan extensions respectively.

In the special case where $f:I^{\star}\to 0^{\star}$ comes from
the unique functor $I\to 0$, we shall write $\holp_I$ for $f_!$
and $\holo_I$ for $f_*$.
}\end{kanax}

\begin{lem}\label{prego}
Let $(f,g)$ be a pair of adjoint functors in $\dia$ and let
   $$\phi:fg\to 1,\ \ \ \psi:1\to gf$$
be the adjunction morphisms. Then $(f^*,g^*)$ is a pair of adjoint
functors and
   $$\phi^*:g^*f^*\to 1,\ \ \ \psi^*:1\to f^*g^*$$
are the adjunction morphisms.
\end{lem}

\begin{proof}
We have the following maps
   $$f\lra{f\psi}fgf\lra{\phi f}f,\ \ \ g\lra{\psi g}gfg\lra{g\phi}g$$
with $(\phi f)(f\psi)=1_f$ and $(g\phi)(\psi g)=1_g$. Then,
   $$(\phi f)^*(f\psi)^*=(f^*\phi^*)(\psi^*f^*)=1_{f^*},\ \ \ (g\phi)^*(\psi g)^*=(\phi^*g^*)(g^*\psi^*)=1_{g^*}$$
whence the assertion.
\end{proof}

\begin{defs}{\rm
We refer to a functor~\eqref{dia} as a {\it left (respectively
right) system of diagram categories\/} if the above axioms (the
Functoriality Axiom, the Isomorphism Axiom, the Disjoint Union
Axiom, and the Left (respectively Right) Homotopy Kan Extension
Axiom) are satisfied. }\end{defs}

In what follows we shall refer to a left and right system of
diagram categories as a {\it bisystem of diagram categories}.

\begin{ex}{\rm
Let $\cc C$ be a closed model category, and let $I\in\dirf$. There
is a natural structure of a closed model category for $\cc C^I$
(see~\cite{F}). Suppose further that $\cc C$ has a zero object.
Denote by $\Ho\cc C^I$ the homotopy category obtained by inverting
the weak equivalences. There is a canonical functor $\cc C^I\to\cc
C^{I^\star}$ which extends a $I$-diagram to $I^\star$ by sending
the zero object and morphisms in $I^\star$ to the zero object and
morphisms in $\cc C$. Any functor $f:I^\star\to J^\star$ therefore
defines a functor $f^*:\cc C^J\to\cc C^I$. It preserves weak
equivalences, hence it defines a functor between homotopy
categories. It follows from~\cite[1.3.2]{F} that the functor
   $$I\in\dirf\longmapsto\Ho\cc C^I$$
determines a bisystem of diagram categories of the domain $\dirf$.

Given an arbitrary model category $\cc C$, let $\cc C_*$ denote
the model category under the terminal object $*$
(see~\cite[p.~4]{Ho}). Then $\cc C_*$ is pointed. Its bisystem of
diagram categories of the domain $\dirf$ is, by definition, that
associated to $\cc C_*$.

}\end{ex}

Let $F:\bb A\to\bb C$ be a morphism between two left systems of
diagram categories $\bb A$ and $\bb C$, and let $f:I\to J$ be a
map in $\dia$. Consider the adjunction maps
   $$\alpha:1\lra{}f^*f_!\ \ \ \textrm{ and }\ \ \ \beta:f_!f^*\lra{}1.$$
Denote by $\gamma_{F,f}$ the composed map
   $$f_!F\xrightarrow{f_!F\alpha}f_!Ff^*f_!\xrightarrow{f_!\iota_{F,f}^{-1}f_!}f_!f^*Ff_!\xrightarrow{\beta Ff_!}Ff_!.$$
We say that $F$ is {\it right exact\/} if $\gamma_{F,f}$ is an isomorphism and if
the following two compatibility relations hold:
   \begin{equation}\label{gr}
    F\alpha_{\bb A}=\iota_{F,f}f_!\circ f^*(\gamma_{F,f})\circ\alpha_{\bb C}F
    \textrm{ and }
    F\beta_{\bb A}=\beta_{\bb C}F\circ f_!(\iota^{-1}_{F,f})\circ\gamma^{-1}_{F,f}f^*.
   \end{equation}
That is, $F\alpha_{\bb A}$ is the composite
   $$F\xrightarrow{\alpha_{\bb C}F}f^*f_!F\xrightarrow{f^*(\gamma_{F,f})}f^*Ff_!\xrightarrow{\iota_{F,f}f_!}Ff^*f_!$$
and $F\beta_{\bb A}$ is the composite
   $$Ff_!f^*\xrightarrow{\gamma^{-1}_{F,f}f^*}f_!Ff^*\xrightarrow{f_!(\iota^{-1}_{F,f})}f_!f^*F\xrightarrow{\beta_{\bb C}F}F.$$
The notion of a left exact (respectively exact) morphism between
two right systems of diagram categories (respectively between two
bisystems of diagram categories) is similarly defined.

\begin{prop}\label{pp}
Let $F:\bb A\to\bb C$ and $G:\bb B\to\bb C$ be two right exact
(respectively left exact) morphisms between left system of diagram
categories (respectively right system of diagram categories); then
the fibred product $\prod(F,G)$ is a left system of diagram
categories (respectively right system of diagram categories) as
well.
\end{prop}

\begin{proof}
We prove the assertion for left systems of diagram categories. The
case of a right system of diagram categories is proved by
symmetry. By Proposition~\ref{qqq} $\prod(F,G)$ is a presystem of
diagram categories. Obviously, it satisfies both the isomorphism
axiom and the disjoint union axiom. We must thus verify the
homotopy Kan extension axioms.

Let $f:I\to J$ be a map in $\dia$. We define the functor
$$f_!:\prod(F,G)_I\to\prod(F,G)_J$$ as follows:
   $$(A,c,B)\longmapsto(f_!(A),\gamma_{G,f}f_!(c)\gamma^{-1}_{F,f},f_!(B)).$$
Then the adjunction maps $\alpha_{\bb A,\bb B}:1\to f^*f_!$ and
$\beta_{\bb A,\bb B}:f_!f^*\to 1$ determine those for $\prod(F,G)$. To see
this we have to check that the squares
   $$\xymatrix{FA\ar[r]^{c}\ar[d]_{F\alpha_{\bb A}} & GB\ar[d]^{G\alpha_{\bb B}}\\
               Ff^*f_!A\ar[r]_{c'}&Gf^*f_!B}$$
with $c'=\iota_{G,f}f_!\circ f^*(\gamma_{G,f})\circ f^*f_!(c)\circ f^*(\gamma^{-1}_{F,f})\circ\iota^{-1}_{F,f}f_!$
and
   $$\xymatrix{Ff_!f^*A\ar[r]^{c''}\ar[d]_{F\beta_{\bb A}} & Gf_!f^*B\ar[d]^{G\beta_{\bb B}}\\
               FA\ar[r]_{c}&GB}$$
with $c''=\gamma_{G,f}f^*\circ f_!(\iota_{G,f})\circ f_!f^*(c)\circ f_!(\iota^{-1}_{F,f})\circ\gamma^{-1}_{F,f}f^*$
are commutative.

We have the following commutative diagram.
   $$\xymatrix{
     FA\ar[r]^{c}\ar[d]_{\alpha_{\bb C}F} & GB\ar[d]^{\alpha_{\bb C}G}\\
     f^*f_!FA\ar[r]^{f^*f_!(c)}\ar[d]_{\iota_{F,f}f_!\circ f^*(\gamma_{F,f})}&f^*f_!GB\ar[d]^{\iota_{G,f}f_!\circ f^*(\gamma_{G,f})}\\
     Ff^*f_!A\ar[r]_{c'}&Gf^*f_!B}$$
In view of the relation~\eqref{gr} one obtains
   \begin{align*}
    G\alpha_{\bb B}\circ c=\iota_{G,f}&f_!\circ f^*(\gamma_{G,f})\circ\alpha_{\bb C}G\circ c=
    \iota_{G,f}f_!\circ f^*(\gamma_{G,f})\circ f^*f_!(c)\circ\alpha_{\bb C}F=\\
    &=\iota_{G,f}f_!\circ f^*(\gamma_{G,f})\circ f^*f_!(c)\circ f^*(\gamma^{-1}_{F,f})\circ
    \iota^{-1}_{F,f}f_!\circ F\alpha_{\bb A}=c'\circ F\alpha_{\bb A}.
   \end{align*}
So, the first square is commutative. Commutativity of the second square is similarly proved.
It is routine to check that both the following composites are the identities (of
$f^*_{\prod(F,G)}$, respectively $f_{_!\prod(F,G)}$).
   $$f^*\lra{\alpha f^*}f^*f_!f^*\lra{f^*\beta}f^*,\ \ \ f_!\lra{f_!\alpha}f_!f^*f_!\lra{\beta f_!}f_!.$$
This yields the homotopy Kan extension axiom.
\end{proof}

\subsection{Consequences of the axioms}

In this section we discuss some consequences of the axioms. We also refer the reader
to Franke's work~\cite{F}.

\subsubsection{Properties of the homotopy Kan extension functors}

A map $f:I\to J$ in $\di$ is a {\it closed (open) immersion\/} if
it is a fully faithful inclusion such that for any $x\in J$ the
relation $\Hom(I,x)\ne\emptyset$ ($\Hom(x,I)\ne\emptyset$) implies
$x\in I$. The following is straightforward.

\begin{immersion}
Let $f:I\to J$ be a closed (respectively open) immersion in $\di$.
Then the map $g:J^\star\to I^\star$ taking $j\in J$ to $j$ if
$j\in I$ and to $\star$ otherwise is a right (respectively left)
adjoint to $f^\star$.
\end{immersion}

\begin{prop}\label{www}
Suppose $\bb C$ is a left system of diagram categories. Let
$f:I\to J$ be a functor, $x\in J$, and let
   \begin{align*}
    i_x:J/x&\to J\\
    j_x:f/x&\to I\\
    l:f/x&\to J/x
   \end{align*}
be the canonical functors. If $J$ is a poset then for $A\in\bb C_I$ we
have isomorphisms
   \begin{align*}
    (f_!A)_x&\iso\holp_{J/x}i^*_xf_!A\\
    {}&\iso\holp_{J/x}l_!j_x^*A\\
    {}&\iso\holp_{f/x}j_x^*A.
   \end{align*}
If $\bb C$ is a right system of diagram categories, a dual
assertion holds for projective homotopy limits and right homotopy
Kan extensions.
\end{prop}

\begin{proof}
Let $c:0\to J$ and $d:0\to J/x$ be the functors taking $0$ to $x$ and $(x,\id)$
respectively. It follows that $d$ is a right adjoint to $J/x\to 0$, and hence
$d^*$ is isomorphic to $\holp_{J/x}$. We then have
   $$f_!A_x=c^*f_!A=d^*i_x^*f_!A\iso\holp_{J/x}i_x^*f_!A$$
whence the first isomorphism follows.

Since $J$ is a poset, $i_x$ is an open immersion and so
$i_x^\star$ has a left adjoint $i_{x+}$ by the immersion lemma. It
sends every object $y\in J$ to $(y,y\le x)$ if $y\le x$, and to
$(\star,\star\le x)$ otherwise. Similarly, $j_x^\star$ has a left
adjoint $j_{x+}$ which sends $y\in I$ to $(\star,\star\le x)$ if
$f(y)\nleqslant x$, and to $(y,f(y)\le x)$ otherwise. It follows
that $i_{x+!}$ and $j_{x+!}$ are isomorphic to $i^*_{x}$ and
$j^*_{x}$ respectively. Since $lj_{x+}=i_{x+}f$, we see that
$l_!j_{x+!}\iso i_{x+!}f_!$, and hence $l_!j^*_{x}$ is isomorphic
to $i^*_{x}f_!$. This implies the second isomorphism. The last
isomorphism is obvious.
\end{proof}

In the case where $f$ is the inclusion of a full subcategory $I\subseteq J$ with $J$
a poset, then for every object $x\in I$ the category $f/x$ has a final object $(x,\id_x)$.
Therefore we have the following

\begin{cor}\label{aaa}
Let $\bb C$ be a left system of diagram categories (respectively
right system of diagram categories), and let $f:I\to J$ be the
inclusion of a full subcategory with $J$ a poset. Then the
canonical morphism in $\bb C_I$
   $$A\to f^*f_!A\textrm{ \ \ (respectively $f^*f_*A\to A$)}$$
is an isomorphism for every object $A$ of $\bb C_I$.
\end{cor}

The last proposition can often be used to reduce assertions about
the functors $f_!$ and $f_*$ to the similar assertions about
$\holp_J$ and $\holo_J$. The following proposition is concerned
with the question of replacing $J$ by a smaller category.

\begin{prop}\label{vvv}
Suppose $\bb C$ is a left system of diagram categories. Let
$i:I^\star\to J^\star$ be some functor (typically the inclusion of
a subcategory). If $i$ has a left adjoint of the form $l^\star$
with $l:J\to I$, then $\holp_JA\iso \holp_Ii^*A$. Dually, suppose
$\bb C$ is a right system of diagram categories. If $i$ has a
right adjoint of the form $r^\star$ for some functor $r:J\to I$,
then $\holo_JA\iso\holo_Ii^*A$.
\end{prop}

\begin{proof}
(1) Let $l$ be a left adjoint of $i$. Then $l_!$ is naturally isomorphic to $i^*$.
It follows that $\holp_J=\holp_I\circ l_!\iso\holp_I\circ i^*$.
\end{proof}

\subsubsection{Cartesian squares} Let $\Box\in\di$
be the poset $\Delta^1\times\Delta^1$, possessing the following elements:
   $$\xymatrix{(0,0) \ar[d] \ar[r] & (0,1) \ar[d] \\
           (1,0) \ar[r]        & (1,1)}$$
where $\to$ stands for $<$.  Let $\la\subset\Box$ be the subposet
obtained by removing the lower right corner $(1,1)$, and let
$\ra\subset\Box$ be the subposet containing all elements of $\Box$
save for $(0,0)$. Let $i_{{}_\ulcorner}:\la\to\Box$ and
$i_{\lrcorner}:\ra\to\Box$ be the inclusions. Let $\bb C$ be a
right system of diagram categories (respectively left system of
diagram categories). An object $A$ of $\bb C_\Box$ is called {\it
cartesian\/} (respectively {\it cocartesian\/}) iff the canonical
morphism $A\to i_{{}_{\lrcorner *}}i^*_{\lrcorner}A$ is an
isomorphism (respectively iff the canonical morphism
$i_{{}_{\ulcorner!}}i_{{}_\ulcorner}^*A\to A$ is an isomorphism).

\begin{lem}\label{hhh}
Let $\bb C$ be a left system of diagram categories. An object $A$
of $\bb C_\Box$ is cocartesian \ifff
$A_{(1,1)}\iso\holp_{{}_{\ulcorner}}i^*_{{}_{\ulcorner}}A$.
Dually, if $\bb C$ is a right system of diagram categories, then
an object $A$ of $\bb C_\Box$ is cartesian \ifff
$A_{(0,0)}\iso\holo_{{\lrcorner}}i^*_{{\lrcorner}}A$.
\end{lem}

\begin{proof}
Let $A$ be an arbitrary object of $\bb C_{\Box}$. By Corollary~\ref{aaa} it follows that
a natural morphism
   $$i^*_{{}_{\ulcorner}}A\to i^*_{{}_{\ulcorner}}i_{{}_{\ulcorner!}}i^*_{{}_{\ulcorner}}A$$
is an isomorphism. Therefore $A_{(i,j)}\iso i_{{}_{\ulcorner!}}i^*_{{}_{\ulcorner}}A_{(i,j)}$
if $(i,j)\in\{(0,0),(0,1),(1,0)\}$.

It follows from Proposition~\ref{www} that
   $$i_{{}_{\ulcorner!}}i^*_{{}_{\ulcorner}}A_{(1,1)}\iso\holp_{{}_{\ulcorner}}i^*_{{}_{\ulcorner}}A.$$
Now our assertion immediately follows.
\end{proof}

From now on all left or right systems of diagram categories are
assumed to be of the domain $\ordd$. Let an object $A$ of $\bb
C_{I\times\Box}$, $I\in\ord$, be called cartesian if
   $$A\to(\id_I\times i_{\lrcorner})_*(\id_I\times i_{\lrcorner})^*A$$
is an isomorphism, and let cocartesiannes be similarly defined,
replacing a right system of diagram categories $\bb C$ by a left
system of diagram categories, $\ra$ by $\la$ and $(\id_I\times
i_{\lrcorner})_*$ by $(\id_I\times i_{{}_\ulcorner})_!$ and
reversing the direction of the arrow. It follows from the
isomorphism axiom and from Proposition~\ref{dgf} below that an
object $A\in\bb C_{I\times\Box}$ is cartesian (respectively
cocartesian) iff the object
$A_{x,\Box}=(i_{x,I}\times\id_\Box)^*A$ is cartesian (respectively
cocartesian) in $\bb C_\Box$ for all $x\in I$.

For any object $I$ of $\ord$, we denote by $\bb C(I)$ the left
system of diagram categories defined as $\bb C(I)_J=\bb C_{I\times
J}$. Here $I$ plays the role of a parameter.

\begin{prop}\label{dgf}
Let $\bb C$ be a left system of diagram categories. Let $f:I\to J$
be a map in $\ordd$. Then $f$ gives the right exact functor
   $$f^*:\bb C(J)\to\bb C(I)$$
induced by $(f\times 1_K)^*:\bb C_{J\times K}\to\bb C_{I\times K}$
with $K\in\ordd$. In particular, the functor respects cocartesian
squares. A dual assertion also holds for right systems of diagram
categories.
\end{prop}

\begin{proof}
Given a map $g:K\to L$ in $\ordd$ we have to show that a natural morphism
   $$\gamma:(1_I\times g)_!(f\times 1_K)^*\to (f\times 1_L)^*(1_J\times g)_!$$
is an isomorphism. Let $A\in\bb C_{J\times K}$ and $(x,y)\in I\times L$.

The map $\phi_{x,y}:1\times g/(x,y)\to g/y$, $((u,u\to x),(v,g(v)\to y))\longmapsto(v,g(v)\to y)$,
has a right adjoint $\psi_{x,y}:g/y\to 1\times g/(x,y)$, $(v,g(v)\to y)\longmapsto((x,x=x),(v,g(v)\to y))$.
It follows from Proposition~\ref{vvv} that $\holp_{1\times g/(x,y)}\iso\holp_{g/y}\psi_{x,y}^*$.

By Proposition~\ref{www} it follows that
   \begin{align*}
    (1_I\times g)_!&(f\times 1_K)^*A_{(x,y)}\iso\holp_{1_I\times g/(x,y)}j^*_{(x,y)}(f\times 1_K)^*A\\
    &\iso\holp_{g/y}\psi_{x,y}^*j^*_{(x,y)}(f\times 1_K)^*A.
   \end{align*}
On the other hand,
   \begin{align*}
    (f\times 1_L)^*&(1_J\times g)_!A_{(x,y)}=i^*_{(f(x),y)}(1_J\times g)_!A
    \iso\holp_{1_J\times g/(f(x),y)}j^*_{(f(x),y)}A\\
    &\iso\holp_{g/y}\psi_{f(x),y}^*j^*_{(f(x),y)}A=\holp_{g/y}\psi_{x,y}^*j^*_{(x,y)}(f\times 1_K)^*A.
   \end{align*}
We have used here the relation $j_{(f(x),y)}\psi_{f(x),y}=(f\times 1_K)j_{(x,y)}\psi_{x,y}$.
So $\gamma$ is an isomorphism. It is routine to
check that the functor of the proposition satisfies the compatibility relations~\eqref{gr}.
\end{proof}

\begin{defs}{\rm
Let $\bb C$ be a right system of diagram categories (respectively
left system of diagram categories). A {\it square\/} in $I$ is a
functor $i:\Box\to I$ which is injective on the set of objects.
Let $A$ be an object of $\bb C_I$; we say that $A$ {\it makes the
square cartesian\/} (respectively {\it cocartesian\/}) if $i^*A$
is cartesian (respectively cocartesian). }\end{defs}

\begin{prop}\label{ccc}
Let $\bb C$ be a left system of diagram categories, and let $i$ be
a square in a poset $I$. If the functor $\la\to(I-i(1,1)/i(1,1))$
possesses a left adjoint and if $A=f_!B$ with $f:J\to I$ a functor
not containing $i(1,1)$ in its image, then $A$ makes $i$
cocartesian. Let $\bb C$ be a right system of diagram categories.
The same holds if the functor $\ra\to(I-i(0,0)\setminus i(0,0))$
possesses a right adjoint and if $A=f_*B$ with $f:J\to I$ a
functor not containing $i(0,0)$ in its image.
\end{prop}

\begin{proof}
It suffices to prove the first assertion. Let $j:I-i(1,1)\to I$ be the inclusion.
By our assumption on the image of $f$, it factors as $J\lra{g}I-i(1,1)\lra{j}I$.
From Corolllary~\ref{aaa} it follows that $j^*A=j^*j_!g_!B\iso g_!B$, whence
$A=j_!g_!B\iso j_!j^*A$. We have
   $$i^*A_{(1,1)}=A_{i(1,1)}\iso\holp_{I-i(1,1)/i(1,1)}h^*A\iso\holp_{{}_{\ulcorner}}i^*_{{}_{\ulcorner}}i^*A,$$
by Proposition~\ref{www} and Proposition~\ref{vvv}(1), where
   $$h:I-i(1,1)/i(1,1)\lra{}I$$
is the canonical functor. Our assertion follows now from Lemma~\ref{hhh}.
\end{proof}

\begin{prop}\label{ddd}
(Concatenation of squares and (co-)cartesiannes) Let $\bb C$ be a
left system of diagram categories (respectively right system of
diagram categories), and let $d_{0,1,2}:\Delta^1\to\Delta^2$ be
the three monotonic injections, and  $A\in\bb
C_{\Delta^1\times\Delta^2}$. Suppose that
$(d_2\times\id_{\Delta^1})^*A\in\bb C_\Box$ is cocartesian
(respectively $(d_0\times\id_{\Delta^1})^*A\in\bb C_\Box$ is
cartesian). Then $(d_0\times\id_{\Delta^1})^*A$ is cocartesian
(respectively $(d_2\times\id_{\Delta^1})^*A\in\bb C_\Box$ is
cartesian) \ifff $(d_1\times\id_{\Delta^1})^*A$ is cocartesian
(respectively $(d_1\times\id_{\Delta^1})^*A\in\bb C_\Box$ is
cartesian).
\end{prop}

\begin{proof}
Let $I=\Delta^1\times\Delta^2-(1,2)$ and $J=I-(1,1)$, and let $j$ and $k$ be the inclusions
of the subposets $I$ and $J$ into $\Delta^1\times\Delta^2$, and let $l:J\to I$ be the
inclusion. It follows from Proposition~\ref{ccc} that $k_!k^*A$ makes both visible squares
in $\Delta^1\times\Delta^2$
   $$\xymatrix{(0,0)\ar[r]\ar[d]&(0,1)\ar[r]\ar[d]&(0,2)\ar[d]\\
               (1,0)\ar[r]&(1,1)\ar[r]&(1,2)}$$
cocartesian. Proposition~\ref{www} then implies the following isomorphism,
   $$k_!k^*A_{(1,1)}\iso\holp_{{}_{\ulcorner}}i^*_{{}_{\ulcorner}}(d_2\times\id_{\Delta^1})^*A.$$
By Proposition~\ref{ccc} we also have
   $$k_!k^*A_{(1,2)}\iso\holp_{{}_{\ulcorner}}i^*_{{}_{\ulcorner}}(d_1\times\id_{\Delta^1})^*A.$$
By our assumption on $(d_2\times\id_{\Delta^1})^*A$, it follows that $j^*A\iso l_!k^*A$
and that $(d_1\times\id_{\Delta^1})^*A$ is cocartesian iff $A$ is isomorphic to
$k_!k^*A=j_!l_!k^*A\iso j_!j^*A$. By Proposition~\ref{ccc} $j_!j^*A$ makes the right square cocartesian.
Since $j_!j^*A_{(1,2)}\iso\holp_{{}_{\ulcorner}}i^*_{{}_{\ulcorner}}(d_0\times\id_{\Delta^1})^*A$
this is the case iff
$(d_0\times\id_{\Delta^1})^*A$ is cocartesian.
\end{proof}

\begin{prop}\label{nm}
Let $\bb C$ be a left system of diagram categories (respectively
right system of diagram categories). For every $I$, the category
$\bb C_I$ has a zero-object and finite coproducts (respectively
products). For every functor $f:I\to J$ the functor $f_!$
(respectively $f_*$) preserves coproducts (respectively products).
\end{prop}

\begin{proof}
Let $f:I^\star\to\emptyset^\star$ be the unique functor. The
inclusion $g:\emptyset^\star\to I^\star$ is left and right adjoint
to $f$. It follows that $g^*$ is left and right adjoint to $f^*$.
Therefore given $0\in\bb C_\emptyset$ the object $f^*0$ (denote it
also by 0) is a zero-object in $\bb C_I$.

Let $I\coprod I$ be the disjoint union of two copies of $I$, and
let $p:I\coprod I\to I$ be the functor which is the identity on
both copies of $I$. By the disjoint union axiom $\bb C_{I\coprod
I}\iso\bb C_I\times\bb C_I$. Hence, the functor $p_!$ provides the
coproduct. Since $f_!$ is left adjoint to $f^*$ with $f:I\to J$ a
map in $\di$, it preserves coproducts.
\end{proof}

Let $f:I^\star\to J^\star$ be a map in $\dia$ and $x\in I$. If $f(x)=\star$ then
$f^*A_x=i^*_{x,I}f^*A=0$ for any $A\in\bb C_J$. Indeed, the composite
$0^\star\lra{i_{x,I}}I^\star\lra{f} J^\star$ factors as
$0^\star\lra{j}\emptyset^\star\lra{l} J^\star$, whence $f^*A_x=j^*l^*A=0$.

\section{D\'erivateurs}

\subsection{Definitions}
Let $\di$ be a category of diagrams. So far we considered only
functors
   $$\bb C:\dia^{\op}\to\CAT$$
evaluated on the category $\dia$. The horizontal morphisms $I\to
J$ in $\dia$ are given by the functors $I^\star\to J^\star$
mapping $\star$ to $\star$. It is also of particular interest to
consider functors
   \begin{equation}\label{we}
    \bb D:\di^{\op}\to\CAT
   \end{equation}
evaluated on the underlying category $\di$. Here we follow the terminology of~\cite{M}.

All the axioms of section~\ref{fgh} can also be reformulated for morphisms
and bimorphisms in $\di$.

\begin{defs}{\rm
A functor~\eqref{we} is called a {\it pred\'erivateur\/} if it
satisfies the Functoriality Axiom. It is a {\it left (respectively
right) d\'erivateur\/} if the Functoriality Axiom, the Isomorphism
Axiom, the Disjoint Union Axiom, the Left (respectively Right)
Homotopy Kan Extension Axiom, and the Left (respectively Right)
Base Change Axiom below are satisfied. }\end{defs}

\begin{baseax}{\rm
Let $f:I\to J$ be a morphism in $\di$ and $x\in J$.
Consider the diagram in $\di$
   $$\xymatrix{\ar @{}[dr] |{\swarrow\alpha_x}
               f/x \ar[d]_{p} \ar[r]^{j_x} & I \ar[d]^{f} \\
               0 \ar[r]_{i_{x,J}}        & J}$$
with $j_x$ a natural map and $\alpha_x$ the bimorphism
   $$fj_x\to i_{x,I}p,\ \ \ \alpha_x:fj_x(y,a:f(y)\to x)=f(y)\lra{a}x=i_{x,J}p(y,a).$$
The bimorphism $\alpha_x$ induces a bimorphism
$\beta_x:p_!j_x^*\to i_{x,I}^*f_!$ which is the composite
   $$p_!j_x^*\to p_!j^*_xf^*f_!\xrightarrow{p_!\alpha_x^*f_!}p_!p^*i^*_xf_!\to i^*_xf_!.$$
The left base change axiom requires $\beta_x$ to be an isomorphism.

By symmetry, the right base change axiom says that the diagram
   $$\xymatrix{\ar @{}[dr] |{\nearrow\gamma_x}
               f\setminus x \ar[d]_{q} \ar[r]^{l_x} & I \ar[d]^{f} \\
               0 \ar[r]_{i_{x,J}}        & J}$$
yields an isomorphism $\delta_x:i^*_{x,I}f_*\to q_*l_x^*$.
}\end{baseax}

We shall refer to a left and right d\'erivateur as a {\it
bid\'erivateur}.

\begin{ex}{\rm
Given a category $\cc C$ closed under colimits, then the
representable pred\'erivateur associated to $\cc C$ is a left
d\'erivateur. A typical example of a bid\'eri\-vateur (of the
domain $\dirf$) is given by the functor
   $$I\longmapsto\Ho\cc C^I$$
with $\cc C$ a closed model category (see~\cite{C} for details).
}\end{ex}

From now on all left or right d\'erivateurs are assumed to be of
the domain $\di$. The notions of a morphism between two
pred\'erivateurs, of the fibred product of a pair of morphisms are
defined similar to those for presystems of diagram categories. It
is similarly proved that the fibred product of two morphisms is a
pred\'erivateur and that it is a left (right) d\'erivateur
whenever both morphisms are right (left) exact.

\begin{prop}\label{tt}
Suppose $\bb D$ is a left d\'erivateur. Let $f:I\to J$ be a
functor in $\di$, $x\in J$. Then for $A\in\bb D_I$ we have an
isomorphism
   $$(f_!A)_x\iso\holp_{f/x}j_x^*A.$$
If $\bb D$ is a right d\'erivateur, a dual assertion holds for
projective homotopy limits and right homotopy Kan extensions.
\end{prop}

\begin{proof}
The assertion immediately follows from the base change axiom. Precisely,
   $$\beta_x^{-1}:(f_!A)_x=i_{x,J}^*f_!A\to\holp_{f/x}j^*_xA$$
is an isomorphism (recall the reader that $\holp_{f/x}=p_!$ by definition).
\end{proof}

In the case where $f$ is the inclusion of a full subcategory $I\subseteq J$, then
for every object $x\in I$ the category $f/x$ has a final object $(x,\id_x)$.
Therefore we have the following

\begin{cor}\label{zz}
Let $\bb D$ be a left d\'erivateur (respectively right
d\'erivateur), and let $f:I\to J$ be the inclusion of a full
subcategory. Then the canonical morphism in $\bb D_I$
   $$A\to f^*f_!A\textrm{ \ \ (respectively $f^*f_*A\to A$)}$$
is an isomorphism for every object $A$ of $\bb D_I$.
\end{cor}

The notion of a (co-)cartesian square is defined as above. Below we formulate
similar statements about (co-)cartesian squares without proofs. They repeat those
in the preceding section word for word.

\begin{lem}\label{hh}
Let $\bb D$ be a left d\'erivateur. An object $A$ of $\bb D_\Box$
is cocartesian \ifff
$A_{(1,1)}\iso\holp_{{}_{\ulcorner}}i^*_{{}_{\ulcorner}}A$.
Dually, if $\bb D$ is a right d\'erivateur, then an object $A$ of
$\bb D_\Box$ is cartesian \ifff
$A_{(0,0)}\iso\holo_{{\lrcorner}}i^*_{{\lrcorner}}A$.
\end{lem}

For any object $I$ of $\di$, we denote by $\bb D(I)$ the left
d\'erivateur defined as $\bb D(I)_J=\bb D_{I\times J}$.

\begin{prop}\label{dgff}
Let $\bb D$ be a left d\'erivateur. Let $f:I\to J$ be a map in
$\di$. Then $f$ gives the right exact functor
   $$f^*:\bb D(J)\to\bb D(I)$$
induced by $(f\times 1_K)^*:\bb D_{J\times K}\to\bb D_{I\times K}$
with $K\in\di$. In particular, the functor respects cocartesian
squares. A dual assertion also holds for right d\'erivateurs.
\end{prop}

\begin{prop}\label{cc}
Let $\bb D$ be a left d\'erivateur, and let $i$ be a square in a
category $I\in\di$. If the functor $\la\to(I-i(1,1)/i(1,1))$
possesses a left adjoint and if $A=f_!B$ with $f:J\to I$ a functor
not containing $i(1,1)$ in its image, then $A$ makes $i$
cocartesian. Let $\bb D$ be a right d\'erivateur. The same holds
if the functor $\ra\to(I-i(0,0)\setminus i(0,0))$ possesses a
right adjoint and if $A=f_*B$ with $f:J\to I$ a functor not
containing $i(0,0)$ in its image.
\end{prop}

\begin{prop}\label{dd}
(Concatenation of squares and (co-)cartesiannes) Let $\bb D$ be a
left d\'erivateur. (respectively right d\'erivateur), and let
$d_{0,1,2}:\Delta^1\to\Delta^2$ be the three monotonic injections,
and  $A\in\bb D_{\Delta^1\times\Delta^2}$. Suppose that
$(d_2\times\id_{\Delta^1})^*A\in\bb D_\Box$ is cocartesian
(respectively $(d_0\times\id_{\Delta^1})^*A\in\bb D_\Box$ is
cartesian). Then $(d_0\times\id_{\Delta^1})^*A$ is cocartesian
(respectively $(d_2\times\id_{\Delta^1})^*A\in\bb D_\Box$ is
cartesian) \ifff $(d_1\times\id_{\Delta^1})^*A$ is cocartesian
(respectively $(d_1\times\id_{\Delta^1})^*A\in\bb D_\Box$ is
cartesian).
\end{prop}

\begin{prop}\label{bnm}
Let $\bb D$ be a left d\'erivateur (respectively right
d\'erivateur). For every $I$, the category $\bb D_I$ has an
initial (respectively final) object and finite coproducts
(respectively products). For every functor $f:I\to J$ the functor
$f_!$ (respectively $f_*$) preserves coproducts (respectively
products).
\end{prop}

\begin{proof}
Let $f:\emptyset\to I$ be the inclusion and $0\in\bb D_\emptyset$.
Since $f_!$ is a left adjoint to $f^*$, it follows that $f_!0$
(denote it also by 0) is an initial object in $\bb D_I$.

Let $I\coprod I$ be the disjoint union of two copies of $I$, and let $p:I\coprod I\to I$
be the functor which is the identity on both copies of $I$. By the disjoint union axiom
$\bb D_{I\coprod I}\iso\bb D_I\times\bb D_I$. Hence, the functor $p_!$ provides the coproduct.
Since $f_!$ is left adjoint to $f^*$ with $f:I\to J$ a map in $\di$, it
preserves coproducts and $f^*$ preserves products (whenever they exist).
\end{proof}

\subsection{Pointed d\'erivateurs}

The d\'erivateurs we work with below must also satisfy some extra
conditions. We start with definitions.

\begin{defs}{\rm
A left d\'erivateur is said to be {\it pointed\/} if it satisfies
the following three conditions:

\begin{enumerate}
\item for any closed immersion $f:I\to J$ in $\di$, the structure
functor $f_!$ possesses a left adjoint $f^?$;

\item for any open immersion $f:I\to J$ in $\di$, the structure functor $f^*$
possesses a right adjoint $f_*$;

\item for any open immersion $f:I\to J$ in $\di$ and any object $x\in J$, the base change
morphism of the diagram
   $$\xymatrix{\ar @{}[dr] |{\nearrow\gamma_x}
               f\setminus x \ar[d]_{q} \ar[r]^{l_x} & I \ar[d]^{f} \\
               0 \ar[r]_{i_{x,J}}        & J}$$
yields an isomorphism $\delta_x:i^*_{x,I}f_*\to q_*l_x^*$.
\end{enumerate}
The corresponding notion for a right d\'erivateur to be pointed is
defined by symmetry. }\end{defs}

We note that for any open immersion $f:I\to J$ in $\di$ and any object $x\in J$
a right adjoint $q_*$ with $q:f\setminus x\to 0$ the unique map always exists. Indeed,
if $x$ is not in $I$ then $f\setminus x=\emptyset$ and $q_*$ exists, because $\emptyset\to 0$
is an open immersion. If $x\in I$ then $f\setminus x$ has an initial object $(x,x=x)$ and
we put $q_*=p^*$ with $0\bl{p}\longmapsto(x,x=x)\in f\setminus x$.

Let $\bb D$ be a left (right) pointed d\'erivateur. Then $\bb D_I$
has a zero object for any $I\in\di$. For the inclusion
$\emptyset\to I$ is both a closed and an open immersion, and
therefore $0=f_!0$ ($0=f_*0$) is also a final (initial) object.
Also, it follows that for every open immersion $f:I\to J$ in $\di$
and any object $x\in J$ the ``value" $f_*A_x$ at $x$, $A\in\bb
D_I$, equals either to $A_x$ if $x\in I$ or to 0
otherwise\label{er}.

In what follows, we refer to a left and right pointed d\'erivateur
as a {\it pointed bid\'erivateur}.

\subsection{Examples}

Given $I\in\dirf$ and a Waldhausen category $\cc C$, the functor
category $\cc C^I$ is a Waldhausen category, too. A map $F\to G$
in $\cc C^I$ is a cofibration (respectively weak equivalence) if
$F(x)\to G(x)$ is so for every $x\in I$.

\begin{defs}{\rm
I. Let $\cc A$ be a category with finite coproducts and an initial
object $e$. Assume that $\cc A$ has two distinguished classes of
maps, called {\it weak equivalences\/} and {\it cofibrations}. A
map is called a {\it trivial cofibration\/} if it is both a weak
equivalence and a cofibration. We call $\cc A$ a {\it category of
cofibrant objects\/} if the following axioms are satisfied.

(A) Let $f$ and $g$ be maps such that $gf$ is defined. If two of
$f$, $g$, $gf$ are weak equivalences then so the third. Any
isomorphism is a weak equivalence.

(B) The composite of two cofibrations is a cofibration. Any
isomorphism is a cofibration.

(C) Given a diagram
   $$A\bl u\longleftarrow C\lra{v} B,$$
with $v$ a cofibration (respectively a trivial cofibration), the
pushout $A\coprod_CB$ exists and the map $A\to A\coprod_CB$ is a
cofibration (respectively a trivial cofibration).

(D) Any map $u$ in $\cc A$ can be factored $u=pi$ with $p$ a weak
equivalence and $i$ a cofibration.

(E) For any object $A$ the map $e\to A$ is a cofibration.

Note that $\cc A$ is a category of cofibrant objects in the sense
of Brown~\cite{B}.

}\end{defs}

For instance, the Waldhausen category of bounded complexes
$C^b(\cc E)$ of an exact category $\cc E$ with weak equivalences
quasi-isomorphisms and cofibrations componentwise admissible
monomorphisms is a category of cofibrant objects.

Let $\cc C$ be a Waldhausen category of cofibrant objects and let
$\Ho\cc C$ denote the category obtained from $\cc C$ by inverting
weak equivalences. One can define the notion of the homotopy for
two maps $f$ and $g$ (see~\cite{B}) and then the category $\pi\cc
C$ with the same objects as $\cc C$ and with $\pi\cc C(A,B)$ by
equal to the quotient of $\cc C(A,B)$ by the equivalence relation
$f\sim g$ defined in terms of the homotopy. Then the class of weak
equivalences in $\pi\cc C$ admits a calculus of left
fractions~\cite{B}. Given $I\in\dirf$, it follows
from~\cite[1.31]{C1} that the functor category $\cc C^I$ is a
Waldhausen category of cofibrant objects.

\begin{thm}[Cisinski~\cite{C1}]\label{cis}
If $\cc C$ is a Waldhausen category of cofibrant objects, then the
hyperfunctor
   $$\bb D\cc C:I\in\dirf\longmapsto\bb D\cc C_I=\Ho\cc C^I$$
determines a left pointed d\'erivateur of the domain $\dirf$.
\end{thm}

\section{The $S.$-construction}\label{s}

Throughout this section $\bb B$ is assumed to be either a left
system of diagram categories (of the domain $\ordd$) or a left
pointed d\'erivateur (of the domain $\di$). Let $\Ar\Delta^n$ be
the poset of pairs $(i,j)$, $0\le i\le j\le n$, where
$(i,j)\le(i',j')$ iff $i\le i'$ and $j\le j'$. Regarded as a
category it may be identified to the category of arrows of
$\Delta^n$.

Given $0\le i<j<k\le n$ let
   \begin{equation}\label{aijk}
    a_{i,j,k}:\Box\to\Ar\Delta^n
   \end{equation}
denote the functor defined as follows:
   $$(0,0)\longmapsto(i,j),\ \ \ (0,1)\longmapsto(i,k),\ \ \ (1,0)\longmapsto(j,j),
     \ \ \ (1,1)\longmapsto(j,k).$$
For any integer $n\geq 0$, we denote by $S_n\bb B$ the
full subcategory of $\bb B_{\Ar\Delta^n}$ that consists of the objects $X$
satisfying the following two conditions:
\begin{itemize}
\item[$\diamond$] for any $i\le n$, $X_{(i,i)}$ is isomorphic to zero in $\bb B_0$;
\item[$\diamond$] for any $0\le i<j<k\le n$, $a_{i,j,k}^*X$ is a cocartesian square if $n>1$.
\end{itemize}
The definition of $S_n\bb B$ is similar to that of $S_n\cc C$,
where $\cc C$ is a Waldhausen category (see~\cite{W} for details).
Note that $S_0\bb B$ is the full subcategory of zero objects in
$\bb B_0$. The category $S_1\bb B$ consists of the objects
$X\in\bb B_{\Delta^2}$ with $X_0$ and $X_2$ isomorphic to zero.

\begin{prop}\label{sss}
Let $n\geq 1$ and let $\ell:\Delta^{n-1}\to\Ar\Delta^n$ be the map taking $j$ to $(0,j+1)$.
Then the functor $\ell^*$ induces an equivalence of categories $S_n\bb B$ and
$\bb B_{\Delta^{n-1}}$.
\end{prop}

\begin{proof}
The proof breaks into two steps.

I. We need some specification both for left systems of diagram
categories and for left pointed d\'erivateurs.

(a) Suppose $\bb B$ is a left system of diagram categories. Let us
consider the following full subcategory $I$ of $\Ar\Delta^n$
   $$\left\{\begin{array}{ccc}
     \xymatrix{(0,1)\ar[d]\ar[r]&(0,2)\ar[r]&\cdots\ar[r]&(0,n)\\ (1,1)}
     \end{array}\right\}  \bigcup_{2\le i\le n}(i,i)$$
as well as the map $g:\Delta^{n-1}\to I$, $j\longmapsto (0,j+1)$.
Since $g$ is an open immersion, it follows from the immersion
lemma that $g$ possesses the left adjoint
$f:I^\star\to\Delta^{n-1\star}$, $(0,j)\longmapsto j-1$ and
$(i,i)\longmapsto\star$. Hence $f^*$ is a right adjoint to $g^*$.

Denote by $\wt{\bb B}_I$ the full subcategory of ${\bb B}_I$
consisting of the objects $X\in\bb B_I$ such that $X_{(i,i)}$,
$1\le i\le n$, are isomorphic to zero. We claim that $f^*$ and
$g^*$ are mutually inverse equivalences between $\bb
B_{\Delta^{n-1}}$ and $\wt{\bb B}_I$. Indeed, $f^*A$ is in
$\wt{\bb B}_I$ for every $A\in\bb B_{\Delta^{n-1}}$ and
$g^*f^*=1$. On the other hand, the adjunction map $B\to f^*g^*B$
is an isomorphism for every $B\in\wt{\bb B}_I$.

(b) Suppose $\bb B$ is a left pointed d\'erivateur. Since $g$ is
an open immersion, then the functor $g^*$ possesses a right
adjoint $g_*$. Let us show that $g^*$ and $g_*$ are mutually
inverse equivalences between $\bb B_{\Delta^{n-1}}$ and $\wt{\bb
B}_I$. Indeed, the adjunction map $g^*g_*\to 1$ is an isomorphism
by Corollary~\ref{zz}. Since $g$ is an open immersion, we see that
$g_*B$ is in $\wt{\bb B}_I$ for all $B\in\bb B_{\Delta^{n-1}}$
(see the corresponding remarks on p.~\pageref{er}). It immediately
follows that the adjunction map $B\to g_*g^*B$ is an isomorphism
for every $B\in\wt{\bb B}_I$.

II. Second, let $h:I\to\Ar\Delta^n$ be the inclusion. It follows from Propositions~\ref{ccc} and~\ref{cc}
that for every $A\in\bb B_I$ the object $h_!A$ makes all visible squares
   $$\xymatrix{(0,0)\ar[r]&(0,1)\ar[d]\ar[r]&(0,2)\ar[d]\ar[r]&\cdots\ar[r]&(0,n-1)\ar[d]\ar[r]&(0,n)\ar[d]\\
                          &(1,1)\ar[r]      &(1,2)\ar[d]\ar[r]&\cdots\ar[r]&(1,n-1)\ar[d]\ar[r]&(1,n)\ar[d]\\
                                           &&(2,2)\ar[r]&\cdots\ar[r]&(2,n-1)\ar[d]\ar[r]&(2,n)\ar[d]\\
                                                                     &&&\ddots&\vdots\ar[d]&\vdots\ar[d]\\
                                                                     &&&&(n-1,n-1)\ar[r]&(n-1,n)\ar[d]\\
                                                                     &&&&&(n,n)}$$
in the category $\Ar\Delta^n$ cocartesian. By Propositions~\ref{ddd} and~\ref{dd} the same is true for all
concatenations of visible squares.

It follows from Propositions~\ref{www} and~\ref{tt} that
$h_!A_{(0,0)}$ is isomorphic to zero. By Corollaries~\ref{aaa}
and~\ref{zz} the canonical morphism
   $$A\to h^*h_!A$$
is an isomorphism for all $A\in\bb B_I$. Let $1\le i\le n$ then for any $A\in\wt{\bb B}_I$
   $$0\iso A_{(i,i)}\iso h^*h_!A_{(i,i)}=h_!A_{(i,i)}.$$
Thus $h_!$ takes an object $A$ of $\wt{\bb B}_I$ to one of
$S_n{\bb B}$. We denote the restriction of $h_!$ to $\wt{\bb B}_I$
by the same letter. To show that
   $$h_!:\wt{\bb B}_I\to S_n\bb B$$
is an equivalence, we must check that the adjunction morphism
   \begin{equation}\label{bv}
    h_!h^*B\to B
   \end{equation}
is an isomorphism for any $B\in S_n\bb B$. By the isomorphism axiom it suffices
to prove that this map is an isomorphism at each point $(i,j)\in\Ar\Delta^n$.
Obviously, it is so at each $(i,j)\in I\cup(0,0)$.

Given $1\le i<j\le n$ let us consider the square
   $$a_{0,i,j}:\Box\to\Ar\Delta^n.$$
We denote the restriction of $a_{0,i,j}$ to $\la$ by $\alpha$. The map~\eqref{bv}
induces the map
   \begin{equation}\label{bb}
    \alpha^*h_!h^*B\to\alpha^*B
   \end{equation}
as well as the map
   \begin{align*}\label{vv}
    h_!h^*B_{(i,j)}\iso\holp_{{}_{\ulcorner}}\alpha^*h_!h^*B\to
    \holp_{{}_{\ulcorner}}\alpha^*B\iso B_{(i,j)}.
   \end{align*}
Here we have used Lemmas~\ref{hhh} and~\ref{hh}. We see that the latter map is an isomorphism whenever
\eqref{bb} is. Since $\im\alpha\subset I$, it follows that~\eqref{bb} is always an isomorphism.
Thus~\eqref{bv} is an isomorphism at each $(i,j)\in\Ar\Delta^n$, and hence $h_!$ and $h^*$
are mutual inverses by the isomorphism axiom.

Since $\ell=hg$ the functor
   $$\ell^*=g^*h^*:S_n\bb B\to\bb B_{\Delta^{n-1}}$$
is an equivalence because both $g^*$ and $h^*$ are so by above. If
$\bb B$ is a left system of diagram categories (respectively a
left pointed d\'erivateur), a quasi-inverse to $\ell^*$ is given
by $h_!f^*$ (respectively by $h_!g_*$). This yields the claim.
\end{proof}

We denote by $\bb B(I)$ the left system of diagram categories or
the left pointed d\'erivateur respectively defined as $\bb
B(I)_J=\bb B_{I\times J}$. Every map $f:I\to J$ yields a functor
$f^*:\bb B(J)\to\bb B(I)$. Below we shall need the following.

\begin{prop}\label{kof}
The structure functor $f^*:\bb B(J)_0=\bb B_J\to\bb B(I)_0=\bb B_I$ respects coproducts.
\end{prop}

\begin{proof}
By definition, the coproduct of two objects $A,B\in\bb B_J$ is the value of the
functor
   $$\bb B_J\times\bb B_J\iso\bb B_{J\coprod J}\lra{p_!}\bb B_J$$
on $(A,B)$, where $p:J\coprod J\to J$ is the canonical map. Propositions~\ref{dgf} and~\ref{dgff}
now imply the assertion.
\end{proof}

By Propositions~\ref{dgf} and~\ref{dgff} $f^*:\bb B(J)\to\bb B(I)$ respects
cocartesian squares. Therefore one obtains a functor (denote it by the same letter)
   $$f^*:S_n\bb B(J)\to S_n\bb B(I),$$
and for any bimorphism $\phi:f\to g$, the bimorphism $\phi^*$
induces a natural transformation of functors
   $$\xymatrix{
     S_n\bb B(J)\ar@/^0.6pc/[r]^{f^*} \ar@/_0.6pc/[r]_{g^*} &S_n\bb B(I)}.$$
We put $\bb S_n\bb B_I=S_n\bb B(I)$. Then $\bb S_n\bb B$ is a
presystem of diagram categories or a pred\'erivateur respectively.
$\bb S_0\bb B$ is trivial and for $n\geq 1$ Proposition~\ref{sss}
implies an equivalence
   $$\bb S_n\bb B\iso\bb B(\Delta^{n-1}).$$
Since $\bb B(\Delta^{n-1})$ is a left system of diagram categories
or a left pointed d\'erivateur respectively, it follows that $\bb
S_n\bb B$ is so as well. We thus obtain a simplicial left system
of diagram categories (respectively a left pointed d\'erivateur)
   $$\bb S.\bb B:\Delta^n\longmapsto \bb S_n\bb B.$$

Consider the following simplicial category:
   $$S.\bb B:\Delta^n\longmapsto S_n\bb B.$$
For any $n\geq 0$, let $iS_n\bb B$ denote the subcategory of $S_n\bb B$ whose
objects are those of $S_n\bb B$ and whose morphisms are isomorphisms in $S_n\bb B$,
and let $i.S_n\bb B$ be the nerve of $iS_n\bb B$. We obtain then
the following bisimplicial object:
   $$i.S.:\Delta^m\times\Delta^n\longmapsto i_mS_n\bb B.$$

\begin{lem}\label{rom}
The space $|i.S.\bb B|$ is connected.
\end{lem}

\begin{proof}
The geometric realization of a bisimplicial set is the diagonal.
If $O_1,O_2\in i_0S_0\bb B$ and $f:O_1\to O_2$ is the unique arrow
in $S_0\bb B$ connecting them, then $A=\sigma_0(f)\in i_1S_1\bb B$
has $\partial_0A=O_2, \partial_1A=O_1$.
\end{proof}

\begin{defs}{\rm
The {\it Grothendieck group\/} $K_0(\bb B)$ is the group generated
by the set of isomorphism classes $[B]$ of objects of $\bb B_0$
with the relations that $[B]=[A]\cdot[C]$ for every $E\in S_2\bb
B$ such that $E_{(0,1)}=A$, $E_{(0,2)}=B$, and $E_{(1,2)}=C$.
}\end{defs}

\begin{lem}\label{444}
$\pi_1|i.S.\bb B|\iso K_0(\bb B)$.
\end{lem}

\begin{proof}
$\pi_1|i.S.\bb B|$ is the free group on $\pi_0|i.S_1\bb B|$ modulo
the relations $d_1(x)=d_2(x)d_0(x)$ for every $x\in\pi_0|i.S_2\bb
B|$. This follows from the fact that $\pi_1|i.S_0\bb B|=0$ and the
Bousfield-Friedlander homotopy spectral sequence~\cite{BF}. We
have that $\pi_0|i.S_1\bb B|$ is the set of isomorphism classes of
objects in $\bb B_0$, $\pi_0|i.S_2\bb B|$ is the set of
isomorphism classes of objects in $S_2\bb B$, and the maps
$d_i:S_2\bb B\to S_1\bb B$ send $E$ to $E_{(1,2)}$, $E_{(0,2)}$
and $E_{(0,1)}$, respectively.
\end{proof}

Let $\cc A$ be an exact category. Its bounded derived category
$D^b(\cc A)$ is constructed as follows (we follow here Keller's
definition~\cite{K}).

Let $H^b(\cc A)$ be the homotopy category of the category of
bounded complexes $\cc C=C^b(\cc A)$, i.e., the quotient category
of $\cc C$ modulo homotopy equivalence. Let $Ac(\cc A)$ denote the
full subcategory of $H^b(\cc A)$ consisting of acyclic complexes.
A complex
   $$X^n\lra{}X^{n+1}\lra{}X^{n+2}$$
is called {\it acyclic\/} if each map $X^n\to X^{n+1}$ decomposes
in $\cc A$ as $X^n\bl{e_n}\twoheadrightarrow
D^n\bl{m_n}\rightarrowtail X^{n+1}$ where $e_n$ is an (admissible)
epimorphism and $m_n$ is an (admissible) monomorphism;
furthermore, $D^n\bl{m_n}\rightarrowtail
X^{n+1}\bl{e_{n+1}}\twoheadrightarrow D^{n+1}$ must be an exact
sequence.

If an exact category is idempotent complete then every
contractible complex is acyclic. Denote by $\cc N=\cc N_{\cc A}$
the full subcategory of $H^b(\cc A)$ whose objects are the
complexes isomorphic in $H^b(\cc A)$ to acyclic complexes. There
is another description of $\cc N$. Let $\cc A\to\tilde{\cc A}$ be
the universal additive functor to an idempotent complete exact
category $\tilde{\cc A}$. It is exact and reflects exactness, and
$\cc A$ is closed under extensions in $\tilde{\cc A}$
(see~\cite[A.9.1]{T}). Then a complex with entries in $\cc A$
belongs to $\cc N$ iff its image in $H^b(\tilde{\cc A})$ is
acyclic. The category $\cc N_{\tilde{\cc A}}=Ac(\tilde{\cc A})$ is
a thick subcategory in $H^b(\tilde{\cc A})$. Note that an object
over $\tilde{\cc A}$ is acyclic iff it has trivial homology
computed in an ambient abelian category. It follows that $\cc N$
is a thick subcategory in $H^b(\cc A)$. Denote by $\Sigma$ the
multiplicative system associated to $\cc N$ and call the elements
of $\Sigma$ {\it quasi-isomorphisms}. A map $s$ is a
quasi-isomorphism iff in any triangle
   $$L\lra{s}M\to N\to L[1]$$
the complex $N$ belongs to $\cc N$.

The derived category is defined as
   $$D^b(\cc A)=H^b(\cc A)/\cc N=H^b(\cc A)[\Sigma^{-1}].$$
Clearly, a map is a quasi-isomorphism iff its image in
$C^b(\tilde{\cc A})$ is a quasi-isomorphism and iff its image in
$D^b(\cc A)$ is an isomorphism.

Recall that the Grothendieck group $K_0(D^b(\cc A))$ is defined as
the group generated by the set of isomorphism classes $[B]$ of
objects of $D^b(\cc A)$ with the relations that $[B]=[A]+[C]$ for
every triangle
   $$A\to B\to C\to A[1].$$

According to~\cite{K} (consult also~\cite{C1}) the hyperfunctor
   $$I\longmapsto D^b(\cc A^I)$$
yields a pointed bid\'erivateur of the domain $\dirf$. It will be
denoted by $\bb D^b(\cc A)$.

\begin{lem}\label{555}
$K_0(\bb D^b(\cc A))=K_0(D^b(\cc A))$.
\end{lem}

\begin{proof}
It is enough to observe that there is a 1-1 correspondence between
the isomorphism classes of objects in $S_2\bb D^b(\cc A)$ and the
isomorphism classes of triangles in $D^b(\cc A)$
(consult~\cite{C1,M1} for details).
\end{proof}

\begin{defs}{\rm
The {\it Algebraic $K$-theory\/} for a small left system of
diagram categories of the domain $\ordd$ or for a left pointed
d\'erivateur $\bb B$ of the domain $\di$ is given by the pointed
space (a fixed zero object 0 of $\bb B_0$ is taken as a basepoint)
   $$K(\bb B)=\Omega|i.S.\bb B|.$$
The $K$-groups of $\bb B$ are the homotopy groups of $K(\bb B)$
   $$K_*(\bb B)=\pi_*(\Omega|i.S.\bb B|)=\pi_{*+1}(|i.S.\bb B|).$$
}\end{defs}

\begin{conv}{\rm
We shall also denote by $0\in\bb B_I$ the object $const^*0$ where
$const:I\to 0$ is the constant map and $0$ is the fixed zero
object of $\bb B_0$. Let $(\textit{L.s.d.c.}$, $\textit{Left\
pointed\ d\'erivateurs})$ denote the corresponding categories of
left systems of diagram categories and left pointed d\'erivateurs
and right exact functors. In order to make the map
   $$(\textit{L.s.d.c.,\ Left\ pointed\ d\'erivateurs})\lra{K}(\textit{Spaces})$$
functorial, in what follows we assume that $\iota_{F,f}:f^*F\to
Ff^*$ are identities for any right exact morphism $F:\bb A\to\bb
B$ and any map $f$ in $\di$.
}\end{conv}

Any right exact functor $F:\bb A\to\bb B$ induces a map $F_*:K(\bb
A)\to K(\bb B)$ of spaces and of their homotopy groups $K_i(\bb
A)\to K_i(\bb B)$.

We can apply the $S.$-construction to each $\bb S_n\bb B$,
obtaining a bisimplicial left system of diagram categories or a
bisimplicial left pointed d\'erivateur respectively. Iterating
this construction, we can form the multisimplicial object $\bb
S.^n\bb B=\bb S.\bb S.\cdots\bb S.\bb B$ and the multisimplicial
categories $iS.^n\bb B$ of isomorphisms. If the additivity theorem
holds, we show that $|i.S.^n\bb B|$ is the loop space of
$|i.S.^{n+1}\bb B|$ for any $n\geq 1$ and that the sequence
   $$\Omega|i.S.\bb B|,\Omega|i.S.S.\bb B|,\ldots,\Omega|i.S.^n\bb B|,\ldots$$
forms a connective $\Omega$-spectrum $\bb K\bb B$. In this case,
one can think of the $K$-theory of $\bb B$ in terms of this
spectrum. This does not affect the $K$-groups, because:
   $$\pi_i(\bb K\bb B)=\pi_i(K(\bb B))=K_i(\bb B),\ \ \ i\geq 0.$$
The additivity theorem remains open for $K(\bb B)$. Nevertheless,
if the definition of $K$-theory as $\Omega|i.S.\bb B|$ is
substituted with the infinite loop space
   $$\Omega^\infty|i.S.^\infty\bb B|= \lim_n\Omega^n|i.S.^n\bb B|,$$
then the additivity theorem does hold (excluding pathological
cases we never have in practice).

Maltsiniotis~\cite{M} uses the $Q.$-construction to define a
$K$-theory space of a triangulated d\'erivateur. This construction
can be extended to arbitrary left systems of diagram categories or
left pointed d\'erivateurs if we replace bicartesian squares in
Maltsiniotis' definition by cocartesian squares. To be more
precise, it is given by the bisimplicial category $Q\bb
B=\{Q_{m,n}\bb B\}_{m,n\geq 0}$ with $Q_{m,n}\bb B$ being the full
subcategory in $\bb B_{\Delta^m\times\Delta^n}$ such that every
$X\in Q_{m,n}\bb B$ makes any square
$i:\Box\to\Delta^m\times\Delta^n$ cocartesian. Let $iQ\bb B$
denote the corresponding maximal groupoid. Then $i.Q\bb B$ is a
trisimplicial object and the $K$-theory space is defined as
$\Omega|\diag(i.Q\bb B)|$.

According to~\cite{C2} the resulting $K$-theory is equivalent to
that defined by us in terms of the $S.$-construction. The proof is
based on~\cite[p.~334]{W} and makes sense without any problem to
our setting.

\section{Simplicial preliminaries}

Multi-simplicial sets will naturally arise in this work. It will
be important that we can work with them directly, without diagonalizing
away all the structure. Such work depends on a couple of lemmas which
we give below. We formulate them for bisimplicial sets as the corresponding lemmas
for multi-simplicial sets are immediate consequences, by taking
suitable diagonals.

\begin{lem}\label{segal}
Let $X..\to Y..$ be a map of bisimplicial sets. Suppose that for every $n$, the
map $X.{}_n\to Y.{}_n$ is a homotopy equivalence. Then $X..\to Y..$ is a homotopy equivalence.
\end{lem}

\begin{proof}
See~\cite{S}.
\end{proof}

\begin{lem}\label{wald}
Let $X..\to Y..\to Z..$ be a sequence of bisimplicial sets so that $X..\to Z..$ is constant.
Suppose that $X.{}_n\to Y.{}_n\to Z.{}_n$ is a fibration up to homotopy, for every $n$.
Suppose further that $Z.{}_n$ is connected for every $n$. Then $X..\to Y..\to Z..$
is a fibration up to homotopy.
\end{lem}

\begin{proof}
\cite[Lemma~5.2]{W1}.
\end{proof}

\begin{lem}\label{lkj}
Let $\cc A$ and $\cc B$ be two small simplicial categories so that
the underlying sets of objects form simplicial sets and let $i\cc
A$ and $i\cc B$ denote the corresponding simplicial subcategories
of isomorphisms. Then every equivalence $F:\cc A\to\cc B$ induces
a homotopy equivalence of bisimplicial objects $F:i.\cc A\to i.\cc
B$. In particular, if $\cc A$ and $\cc B$ happen to be two left
systems of diagram categories or two left pointed d\'erivateurs,
then every right exact equivalence $F:\cc A\to\cc B$ induces a
homotopy equivalence $F:i.S.\cc A\to i.S.\cc B$.
\end{lem}

\begin{proof}
By Lemma~\ref{segal} it suffices to show that each $F_n:i.\cc A_n\to i.\cc B_n$ is a
homotopy equivalence. The latter is obvious, because $F_n:i\cc A_n\to i\cc B_n$ is
an equivalence of categories.

If $F:\cc A\to\cc B$ is an equivalence of left systems of diagram
categories or left pointed d\'erivateurs, it is directly verified
that $F$ induces an equivalence of categories $F:S_n\cc A\to
S_n\cc B$, and hence a homotopy equivalence $F:i.S.\cc A\to
i.S.\cc B$.
\end{proof}

Let $C$ and $D$ be two simplicial objects in a category $\cc C$
and let $\Delta/\Delta^1$ denote the category of objects over $\Delta^1$ in
$\Delta$; the objects are the maps $\Delta^n\lra{}\Delta^1$. For any
simplicial object $C$ in $\cc C$ let $C^*$ denote the composed
functor
   \begin{gather*}
    (\Delta/\Delta^1)^{\op}\lra{}\Delta^{\op}\lra{C}\cc C\\
    (\Delta^n\lra{}\Delta^1)\longmapsto\Delta^n\longmapsto C_n.
   \end{gather*}
Then a {\it simplicial homotopy\/} of maps from $C$ to $D$ is a
natural transformation $C^*\lra{}D^*$~\cite[p.~335]{W}.

There is a functor $P:\Delta\to\Delta$ with
$P\Delta^n=\Delta^{n+1}$ such that the natural map
$s_0:\Delta^n\to\Delta^{n+1}=P\Delta^n$ is a natural
transformation $\id_\Delta\to P$. It is obtained by formally
adding an initial object $0'$ to each $\Delta^n$ and then
identifying $\{0'<1<\cdots<n\}$ with $\Delta^{n+1}$. Thus
$P(s_i)=s_{i+1}$ and $P(d_i)=d_{i+1}$. If $A$ is a simplicial
object in $\cc A$, the {\it path space\/} $PA$ is the simplicial
object obtained by composing $A$ with $P$. Thus $PA_n=A_{n+1}$,
and the $i$th face operator on $PA$ is the $\partial_{i+1}$ of
$A$, and the $i$th degeneracy operator on $PA$ is the
$\sigma_{i+1}$ of $A$. Moreover, the maps $\partial_0:A_{n+1}\to
A_n$ form a simplicial map $PA\to A$.

Let us write $A_0$ for the constant simplicial object at $A_0$. The natural
maps $\sigma_0^{n+1}:A_0\to A_{n+1}$ form a simplicial map $\iota:A_0\to PA$,
and the maps $A_{n+1}\to A_0$ induced by the canonical inclusion of $\Delta^0=\{0\}$
in $\Delta^{n+1}$ form a simplicial map $\rho:PA\to A_0$ such that $\rho\iota$
is the identity on $A_0$. On the other hand, $\iota\rho$ is homotopic to
the identity on $PA$. The homotopy is given by the natural transformation
   $$(a:\Delta^n\to\Delta^1)\longmapsto(\phi_a^*:A_{n+1}\to A_{n+1})$$
induced from $(a:\Delta^n\to\Delta^1)\longmapsto(\phi_a:\Delta^{n+1}\to\Delta^{n+1})$
where $\phi_a(0)=0$ and
   $$\phi_a(j+1)=
     \left\{
      \begin{array}{rcl}
       j+1,\ a(j)&=&1\\
       0,\ a(j)&=&0.
      \end{array}
      \right.$$
This shows that $PA$ is homotopy equivalent to $A_0$.

\section{$\Gamma$-spaces}

In this section we use Segal's machine~\cite{S} to get some
information about the $K$-theory $K(\bb B)$. We start with
preparations.

Given a finite set $T$ by $\cc P(T)$ denote the set of subsets of
$T$ and the set $\{1,2,\ldots,n\}$ is denoted by $\bb n$.

\begin{defs}{\rm
I. $\Gamma$ is the category whose objects are all finite sets, and
whose morphisms from $S$ to $T$ are the maps $\theta:S\to\cc P(T)$
such that $\theta(\alpha)$ and $\theta(\beta)$ are disjoint when
$\alpha\ne\beta$. The composite of $\theta:S\to\cc P(T)$ and
$\phi:T\to\cc P(U)$ is $\psi:S\to\cc P(U)$, where
$\psi(\alpha)=\bigcup_{\beta\in\theta(\alpha)}\phi(\beta)$.

II. A {$\Gamma$-space\/} is a contravariant functor $A$ from
$\Gamma$ to spaces such that

(a) $A(\bb 0)$ is contractible, and

(b) for any $n$ the map $p_n:A(\bb n)\to A(\bb
1)\times\bl{n}\cdots\times A(\bb 1)$ induced by the maps $i_k:\bb
1\to\bb n$ in $\Gamma$, where $i_k(1)=\{k\}\subset\bb n$, is a
homotopy equivalence.

We shall refer to $A(\bb 1)$ as the {\it underlying space}.

}\end{defs}

There is a covariant functor $\Delta\to\Gamma$ which takes
$\Delta^m$ to $\bb m$ and $f:\Delta^m\to\Delta^n$ to
$\theta(i)=\{j\in\bb n\mid f(i-1)<j\le f(i)\}$. Using this functor
one can regard $\Gamma$-spaces as simplicial spaces.

Segal uses a realization functor $A\to|A|$ for simplicial spaces
which is slightly different from the usual one
(see~\cite[Appendix~A]{S}). If $A$ is a $\Gamma$-space its
realization will mean the realization of the simplicial space it
defines.

\begin{defs}{\rm
If $A$ is a $\Gamma$-space, its classifying space is the
$\Gamma$-space $BA$ such that, for any finite set $S$, $BA(S)$ is
the realization of the $\Gamma$-space $T\longmapsto A(S\times T)$.
}\end{defs}

If $A$ is a $\Gamma$-space the spaces $A(\bb 1),BA(\bb 1),B^2A(\bb
1),\ldots$ form a spectrum, denoted by $\bb BA$. The reason of
introducing $\Gamma$-spaces is that they arise naturally from
categories.

\begin{defs}{\rm
A {$\Gamma$-category\/} is a contravariant functor $\cc C$ from
$\Gamma$ to categories such that

(a) $\cc C(\bb 0)$ is equivalent to the category with one object
and one morphism;

(b) for each $n$ the functor $p_n:\cc C(\bb n)\to\cc C(\bb
1)\times\bl{n}\cdots\times\cc C(\bb 1)$ induced by the maps
$i_k:\bb 1\to\bb n$ in $\Gamma$ is an equivalence of categories.
}\end{defs}

If $\cc C$ is a $\Gamma$-category, $|\cc C|$ is a $\Gamma$-space.
Here $|\cc C|$ means the functor $S\longmapsto|\cc C(S)|$.

$\Gamma$-categories arise in the following way. Let $\cc C$ be a
category with a zero object 0 in which sums exist. If $S$ is a
finite set, let $\cc P(S)$ denote the category of subsets of $S$
and their inclusions~--- this should not cause confusion with the
earlier use of $\cc P(S)$. Let $\cc C(S)$ denote the category
whose objects are the functors from $\cc P(S)$ to which take
disjoint unions to sums. For example, $\cc C(\bb 3)$ is the
category of diagrams in $\cc C$ of the form
  $$
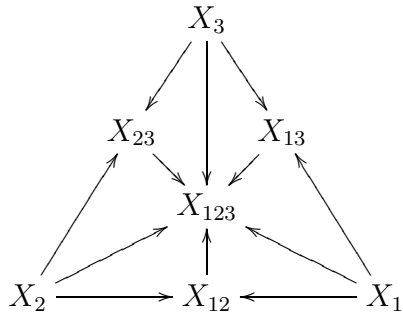
\xy
    *+{X_{23}}="X_{23}",<1cm,1.5cm>*+{X_3}="X_3",-<-1cm,1.5cm>*+{X_{13}}="X_{13}",
    -<1cm,1cm>*+{X_{123}}="X_{123}",
    -<0cm,1.2cm>*+{X_{12}}="X_{12}",-<-2.37cm,0cm>*+{X_1}="X_1",-<4.75cm,0cm>*+{X_2}="X_2",
    \ar "X_3";"X_{23}",\ar "X_3";"X_{13}",\ar "X_3";"X_{123}",\ar "X_{23}";"X_{123}",\ar "X_{13}";"X_{123}",
    \ar "X_2";"X_{12}",\ar "X_1";"X_{12}",\ar "X_2";"X_{23}",\ar "X_{12}";"X_{123}",\ar "X_1";"X_{13}",
    \ar "X_2";"X_{123}",\ar "X_1";"X_{123}"
   \endxy$$
in which each straight line (such as $X_1\to X_{123}\longleftarrow
X_{23}$) is an expression of the middle object as a sum of ends.
The morphisms of $\Gamma$ were so defined that the morphisms from
$S$ to $T$ in $\Gamma$ correspond precisely to functors from $\cc
P(S)$ to $\cc P(T)$ which preserve disjoint unions. The category
$i\cc C(S)$ of isomorphisms in $\cc C(S)$ satisfies the definition
above because, for example, the forgetful functor $i\cc C(\bb
3)\to i\cc C\times i\cc C\times i\cc C$, which takes the diagram
above to $(X_1,X_2,X_3)$ is an equivalence of categories.

Conforming to the terminology and notation
of~\cite[section~1.8]{W}, we denote the resulting simplicial
category by $N.\cc C$. By definition, $N_0\cc C=0$ and $N_n\cc
C=\cc C(\bb n)$ for each $n\geq 1$. We refer to the simplicial
category $N.\cc C$ as the {\it nerve with respect to the
composition law}. By construction, the space $|i.N.\cc C|$ is
$B|i.\cc C|(\bb 1)$. By $N.\bb B$ will be denoted the nerve with
respect to the composition law associated to the category $\bb
B_0$.

Observe that any functor $f:\cc C\to\cc D$ respecting sums yields
a map of bisimplicial objects
   $$f^*:i.N.\cc C\to i.N.\cc D.$$
It follows that for a given left system of diagram categories or a
left pointed d\'erivateur $\bb B$ one can also produce the
multisimplicial categories $iN.^mS.^n\bb B$, $m,n\geq 0$, and the
spaces $|i.N.^mS.^n\bb B|$ by iterating the $N.$- and
$S.$-constructions.

\begin{prop}\label{cosi}
$|i.S.\bb B|$ is canonically an infinite loop space, and hence is
so the $K$-theory space $K(\bb B)$.
\end{prop}

\begin{proof}
The above considerations show that $|i.S.\bb B|$ is the underlying
space of a $\Gamma$-space, with respect to the composition law
produced by coproduct.
\end{proof}

Let $F:\bb A\to\bb B$ be a right exact functor between left
systems of diagram categories or left pointed d\'erivateurs,
respectively. Let further $N_n(\bb A\to\bb B)$ denote the fibred
product of the diagram
   $$N_n\bb A\lra{F}N_n\bb B\bl{\partial_0}\longleftarrow(PN.\bb B)_n=N_{n+1}\bb B.$$
An object of $N_n(\bb A\to\bb B)$ is a triple $(A,c,B)$ with $A\in
N_n\bb A$, $B\in N_{n+1}\bb B$, $c:F(A)\to\partial_0(B)$ an
isomorphism in $N_n\bb B$. One obtains a simplicial category
   $$N.(\bb A\to\bb B):\Delta^n\longmapsto N_n(\bb A\to\bb B).$$
For every $n$, there is a functor
   $$g:\bb B_0=N_1\bb B\to N_n(\bb A\to\bb B)$$
defined by $B\longmapsto(0,1,v^*B)$ with
$v:\Delta^{n+1}\to\Delta^1$, $i\longmapsto 0$ if $i=0$ and
$i\longmapsto 1$ otherwise.

Regarding $\bb B_0$ as a trivial simplicial category, we obtain a
sequence
   $$\bb B_0\lra{g}N.(\bb A\to\bb B)\lra{p}N.\bb A$$
where $p$ is the projection. The latter sequence yields the
sequence
   \begin{equation}\label{si}
    i.S.\bb B\lra{g}i.N.S.(\bb A\to\bb B)\lra{p} i.N.S.\bb A
   \end{equation}
with $N.S.(\bb A\to\bb B)=N.(S.\bb A\to S.\bb B)$. Note that the
space $|i.N.S.\bb A|$ is $B|i.S.\bb A|(\bb 1)$, where $B|i.S.\bb
A|$ is the $\Gamma$-space associated to $|i.S.\bb A|$.

\begin{lem}\label{bene}
The sequence~\eqref{si} is a fibration up to homotopy.
\end{lem}

\begin{proof}
By Lemma~\ref{wald} it is enough to show that for every $n$ the
sequence $i.S.\bb B\lra{}i.N_nS.(\bb A\to\bb B)\lra{}i.N_nS.\bb A$
is a fibration (since the base term $i.N_nS.\bb A=i.\Hom(\cc P(\bb
n),S.\bb A)\iso (i.S.\bb A)^n$ is connected for every $n$ by
Lemma~\ref{rom}). We will show that the sequence is the same, up
to homotopy, as the trivial fibration sequence associated to the
product $i.S.\bb B\times i.N_nS.\bb A$.

Let $u:\Delta^1\to\Delta^{n+1}$ be the map $0;1\longmapsto 0;1$.
Also, consider the maps $d_0:\Delta^n\to\Delta^{n+1}$ and
$s_0:\Delta^{n+1}\to\Delta^n$. We construct the following diagram
for any $B\in N_{n+1}\bb B$,
   $$B'=v^*u^*B\lra{\phi}B\bl{\psi}\longleftarrow B''=s_0^*\partial_0B.$$
For any subset $S$ of $[\bb{n+1}]$,
   $$B'_S=
     \left\{
      \begin{array}{rcl}
       B_1,\ 1&\in&S\\
       0,\ 1&\notin&S
      \end{array}
      \right.$$
and
   $$B''_S=
     \left\{
      \begin{array}{rcl}
       B_{S\setminus\{1\}},\ 1&\in&S\\
       B_S,\ 1&\notin&S
      \end{array}
      \right.$$
whence the definitions of $\phi$ and $\psi$ follow. Note that
$B'_S\lra{\phi_S}B_S\bl{\psi_S}\longleftarrow B''_S$ belongs to
$N_2\bb B$.



The map $N_n(\bb A\to\bb B)\to N_n\bb A\times\bb B_0$,
$(A,c,B)\longmapsto(A,B_{\{1\}})$, is an equivalence of
categories. A quasi-inverse is given by the functor
   $$(A,B)\longmapsto(A,1,s_0^*FA\coprod v^*B).$$
Thus the induced map $i.N_nS.(\bb A\to\bb B)\to i.N_nS.\bb A\times
i.S.\bb B$ is a homotopy equivalence by Lemma~\ref{segal}.

This homotopy equivalence fits into the following commutative
diagram
   $$\begin{CD}
      i.S.\bb B@>>>i.N_nS.(\bb A\to\bb B)@>>>i.N_nS.\bb A\\
      @V1VV@VVV@VV1V\\
      i.S.\bb B@>>>i.N_nS.\bb A\times i.S.\bb B@>>>i.N_nS.\bb A
     \end{CD}$$
Being homotopy equivalent to the trivial fibration (the lower row
of the diagram), we conclude that the upper sequence is a
fibration, as was to be shown.
\end{proof}

As above, one can construct the sequence
   $$i.\bb B_0\to P(i.N.\bb B)\to i.N.\bb B.$$
The composite map is constant and the middle term is contractible,
so we obtain a map well defined up to homotopy,
   $$|i.\bb B_0|\to\Omega|i.N.\bb B|.$$
By naturality we can substitute $\bb B$ with the simplicial
category $S.\bb B$ in the above sequence. We obtain a sequence
   $$i.S.\bb B\to P(i.N.S.\bb B)\to i.N.S.\bb B$$
where the ``$P$" refers to the $N.$-direction. It follows from the
preceding lemma that the sequence is a fibration up to homotopy.
Thus $|i.S.\bb B|\to\Omega|i.N.S.\bb B|$ is a homotopy equivalence
and more generally therefore, in view of Lemma~\ref{segal}, also
the map $|i.N.^nS.\bb B|\to\Omega|i.N.^{n+1}S.\bb B|$. There
results a spectrum
   $$n\longmapsto|i.N.^nS.\bb B|,$$
which is actually a $\Omega$-spectrum. It is nothing but the
spectrum $n\longmapsto B^n|i.S.\bb B|(\bb 1)$ produced by Segal's
machine.

As in the sequence
   $$|i.S.\bb B|\to\Omega|i.N.S.\bb B|\to\Omega\Omega|i.N.N.S.\bb B|\to\cdots$$
all the maps are homotopy equivalences, then so is the map
   $$|i.S.\bb B|\to\Omega^\infty|i.N.^\infty S.\bb
     B|=\lim_n\Omega^n|i.N.^n S.\bb B|$$

\begin{cor}\label{1111}
Suppose we are given a sequence $\bb A\to\bb B\to\bb C$ of right
exact morphisms of left systems of diagram categories or left
pointed d\'erivateurs respectively. Then the square
   $$\begin{CD}
      i.S.\bb B@>>>i.N.S.(\bb A\to\bb B)\\
      @VVV@VVV\\
      i.S.\bb C@>>>i.N.S.(\bb A\to\bb C)
     \end{CD}$$
is homotopy cartesian.
\end{cor}

\begin{proof}
There is a commutative diagram
   $$\begin{CD}
      i.S.\bb B@>>>i.N.S.(\bb A\to\bb B)@>>>i.N.S.\bb A\\
      @VVV@VVV@VV{\id}V\\
      i.S.\bb C@>>>i.N.S.(\bb A\to\bb C)@>>>i.N.S.\bb A
     \end{CD}$$
in which the rows are fibrations up to homotopy by
Lemma~\ref{bene}. Therefore the square on the left is homotopy
cartesian.
\end{proof}

\begin{cor}\label{3333}
The following two assertions are valid.

$(1)$ To a right exact morphism there is associated a fibration
   $$i.S.\bb B\to i.S.\bb C\to i.N.S.(\bb B\to\bb C).$$

$(2)$ If $\bb C$ is a retract of $\bb B$ (by right exact functors)
there is a splitting
   $$i.S.\bb B\iso i.S.\bb C\times i.N.S.(\bb C\to\bb B).$$
\end{cor}

\begin{proof}
(1). If $\bb A=\bb B$ the space $|i.N.S.(\bb A=\bb A)|$ is
contractible whence the first assertion.

(2). This is the case of Corollary~\ref{1111} where the composed
map $\bb A\to\bb B\to\bb C$ is an identity map since $i.N.S.(\bb
A\to\bb C)$ is contractible in that case.
\end{proof}

\section{The additivity theorem}

Let $\bb B$ be either a left system of diagram categories or a
left pointed d\'erivateur. Denote by $\bb E_0$ the full
subcategory in $\bb B_\Box$ consisting of the cocartesian squares
$B\in\bb B_\Box$ with $B_{(1,0)}$ isomorphic to zero. If we
replace $\bb B$ by $\bb B(I)$, we define the category $\bb E_I$
similar to $\bb E_0$. One obtains a left system of diagram
categories or a left pointed d\'erivateur $\bb E$ respectively.

\begin{lem}\label{cxc}
The map $l:\Delta^1\to\Box$, $i\longmapsto(0,i)$, induces an
equivalence of categories $l^*:\bb E_0\to\bb B_{\Delta^1}$. It
also induces a right exact equivalence $\bb E\to\bb B(\Delta^1)$.
\end{lem}

\begin{proof}
The map $l$ factors as
   $$\Delta^1\lra{g}\la\lra{h}\Box.$$
The proof of Proposition~\ref{sss} shows that $l^*:\bb E_0\to\bb
B_{\Delta^1}$ is an equivalence (for left pointed d\'erivateurs
use the fact that $g$ is an open immersion). Obviously, the
induced morphism $\bb E\to\bb B(\Delta^1)$ is an equivalence. It
is right exact by Propositions~\ref{dgf} and~\ref{dgff}.
\end{proof}

\begin{cor}\label{cx}
The map $f:\Box\to\Ar\Delta^2$, $(i,j)\longmapsto(i,j+1)$, induces
an equivalence of categories $f^*:S_2\bb B\to\bb E_0$. It also
induces a right exact equivalence $\bb S_2\bb B\lra{}\bb E$.
\end{cor}

\begin{proof}
Let $\ell:\Delta^1\to\Ar\Delta^2$ be the map $i\longmapsto (0,i+1)$. It factors as
   $$\Delta^1\lra{l}\Box\lra{f}\Ar\Delta^2$$
where $l$ is the map of Lemma~\ref{cxc}. By Proposition~\ref{sss} it follows that
$\ell^*=l^*f^*:\bb S_2\bb B\to\bb B_{\Delta^1}$ is an equivalence. By Lemma~\ref{cxc} $l^*$ is an equivalence,
and hence $f^*$ is so.
\end{proof}

Below we shall need the following.

\begin{lem}\label{iop}
Let $\bb B$ be either a left system of diagram categories or a
left pointed d\'erivateur, and let $B\in\bb B_\Box$ be a
cocartesian square such that the map $B_{(0,0)}\to B_{(0,1)}$
(respectively the map $B_{(0,0)}\to B_{(1,0)}$) is an isomorphism.
Then $B_{(1,0)}\to B_{(1,1)}$ (respectively $B_{(0,1)}\to
B_{(1,1)}$) is an isomorphism as well. On the other hand, a square
with two parallel arrows being isomorphisms is cocartesian.
\end{lem}

\begin{proof}
Suppose that the map $B_{(0,0)}\to B_{(0,1)}$ is an isomorphism.
Let $q:\Box\to\Delta^1$ denote the functor $(\epsilon,\eta)\longmapsto\epsilon$, and
let $i:\Delta^1\to\Box$ be the functor $\kappa\longmapsto(\kappa,0)$. Then $i$ is a left adjoint
to $q$ and hence $i_!\iso q^*$. The map $i$ factors as
   $$\Delta^1\lra{l}\la\lra{i_{{}_\ulcorner}}\Box$$
where $l(\kappa)=(\kappa,0)$. By Propositions~\ref{ccc}
and~\ref{cc} the object $i_!i^*B\iso i_{{}_\ulcorner !}(l_!i^*B)$
is cocartesian.

Let $\beta:(i_!i^* B\iso)q^*i^*B\to B$ be the adjunction morphism. Then $\beta_{(0,0)}=\beta_{(1,0)}=1$ and
$\beta_{(0,1)}$ is an isomorphism by assumption. It follows that $i_{{}_\ulcorner}^*\beta$
is an isomorphism. Consider the following commutative square:
   $$\xymatrix{i_{{}_\ulcorner!}i_{{}_\ulcorner}^*i_!i^*B\ar[r]^{i_{{}_\ulcorner!}i_{{}_\ulcorner}^*\beta}\ar[d]
               &i_{{}_\ulcorner!}i_{{}_\ulcorner}^*B\ar[d]\\
               i_!i^*B\ar[r]^{\beta}& B}$$
The upper arrow is an isomorphism. The vertical maps are isomorphisms too, because
both $B$ and $i_!i^*B$ are cocartesian. We see that $\beta$ is an isomorphism as well.
This implies that $B_{(1,0)}\to B_{(1,1)}\iso\beta_{(1,1)}$ is an isomorphism.
The corresponding assertion when $B_{(0,0)}\to B_{(1,0)}$ is an isomorphism is deduced
from the first one by applying the autoequivalence $\tau:\Box\to\Box$ transposing
the vertices $(1,0)$ and $(0,1)$.

On the other hand, if both $B_{(0,0)}\to B_{(0,1)}$ and $B_{(1,0)}\to B_{(1,1)}$ are
isomorphisms, it follows that $\beta:i_!i^*B\to B$ is an isomorphism. Since $i_!i^*B$
is cocartesian, it follows that $B$ is cocartesian, too.
\end{proof}

We want to construct a map $\alpha:\bb B_0\to\bb E_0$ that takes an object $A\in\bb B_0$
to one in $\bb E_0$ which is depicted in $\bb B_0$ as the square
   $$\xymatrix{A\ar[r]^1\ar[d]&A\ar[d]\\
               0\ar[r]&0.}$$
Let $s_0:\Delta^1\to 0$ be the unique map. First, suppose that
$\bb B$ is a left system of diagram categories and
$d:\Box^\star\to\Delta^{1\star}$ is the map $(0,i)\longmapsto i$
and $(1,i)\longmapsto\star$. We put then $\alpha=d^*s_0^*$. If
$\bb B$ is a left pointed d\'erivateur, let $l:\Delta^1\to\Box$ be
the map $i\longmapsto (0,i)$. Then $l$ is an open immersion, and
hence there is a right adjoint functor $l_*$ to $l^*$. In this
case $\alpha=l_*s_0^*$.

Let $j:\Box\to\Delta^1$ be the map
$(\epsilon,\eta)\longmapsto\eta$. The morphism $j^*$ takes
$B\in\bb B_{\Delta^1}$ to a square in $\bb B_\Box$ which is
depicted in $\bb B_0$ as
   $$\xymatrix{B_0\ar[r]\ar[d]_1&B_1\ar[d]^1\\
               B_0\ar[r]&B_1.}$$

Let $\bb B$ be a left system of diagram categories and let
$u:\Delta^{1\star}\to 0^\star$ be the map $0\longmapsto\star$ and
$1\longmapsto 0$. Then $u^*$ takes an object $B\in\bb B_0$ to one
in $\bb B_{\Delta^1}$ with $u^*B_0=0$ and $u^*B_1=B$. We put
$\beta=j^*u^*$. In turn, if $\bb B$ is a left pointed
d\'erivateur, consider the map $v:0\to\Delta^1$ with $v(0)=1$.
Then $v_!B_0=0$ and $v_!B_1=B$ for any $B\in\bb B_0$. In this case
$\beta:=j^*v_!$.

The $\beta$ takes an object $B\in\bb B_0$ to the square
   $$\xymatrix{0\ar[r]\ar[d]&B\ar[d]^1\\
               0\ar[r]&B}$$

Let $\bb B$ be either a left system of diagram categories or a
left pointed d\'erivateur, and let $\bb B'$ and $\bb B''$ be
either two left subsystems of diagram categories or two left
pointed subd\'erivateurs respectively in such a way that the
inclusion morphisms are right exact. There are three natural right
exact morphisms $s,t,q:\bb E\to\bb B$ taking an object $E\in\bb
E_?$ to $E_{(0,0)}$, $E_{(0,1)}$, and $E_{(1,1)}$ respectively.
Define $\bb E(\bb B',\bb B,\bb B'')$ as a presystem of diagram
categories or as a pred\'erivateur respectively that consists of
the squares $E\in\bb E$ with $E_{(0,0)}\in\bb B'$ and
$E_{(1,1)}\in\bb B''$. Then $\bb E(\bb B',\bb B,\bb B'')$ is a
left system of diagram categories or a left pointed d\'erivateur
respectively. For $\bb E(\bb B',\bb B,\bb B'')$ is equivalent to
the fibred product of the diagram $\bb E\xrightarrow{(s,q)}\bb
B\times\bb B\longleftarrow\bb B'\times\bb B''$. We note that $\bb
E=\bb E(\bb B,\bb B,\bb B)$.

In order to have the unitial and associative $H$-space structure
to $|i.S.\bb B|$ induced by coproduct $\coprod$ via the map
   $$|i.S.\bb B|\times|i.S.\bb B|\lra{\sim}|i.S.\bb B\times i.S.\bb B|\lra{\coprod}|i.S.\bb B|,$$
we must have good choices for $A\coprod B$, $A,B\in S_n\bb B$, in
such a way that $f^*(A\coprod B)=f^*(A)\coprod f^*(B)$ where
$f:\Delta^m\to\Delta^n$ is a structure map in $\Delta$ (one always
has an isomorphism between them because $f^*$ respects coproducts
by Lemma~\ref{kof}). We would have then a simplicial equivalence
$\coprod(\coprod\times 1)\iso\coprod(1\times\coprod)$
   $$\xymatrix{i.S.\bb B\times i.S.\bb B\times i.S.\bb B\ar[rr]^{\coprod\times 1}\ar[dd]_{1\times\coprod}
               && i.S.\bb B\times i.S.\bb B\ar[d]^{\coprod} \\
               &  &i.S.\bb B\ar@/^/@2{->}[dl]^{\alpha}\\
               i.S.\bb B\times i.S.\bb B\ar[r]_{\coprod}&i.S.\bb B}$$
inducing a homotopy between them after realization. It would also
follow that the two maps $i.S.\bb B\to i.S.\bb B$ given by
$B\longmapsto B\coprod 0$ and $B\longmapsto 0\coprod B$ are
homotopic to the identity map, hence $|i.S.\bb B|$ is unitial. It
seems that we do not have enough data to produce such choices in
general. We shall refer to this case as {\it pathological}. The
latter term is caused by the observation that one always has the
required choices in practice. Indeed, all left systems of diagram
categories or left pointed d\'erivateurs arise in practice as the
hyperfunctor
   $$I\longmapsto\Ho\cc C^I$$
with $\cc C$ being closed under coproducts. Then the choices are
made in $\cc C$.

\begin{conv}{\rm
In the rest of this section we assume $\bb B$ to be
non-pathological.
}\end{conv}

By a {\it right exact sequence\/} $F'\to F\to F''$ of right exact
functors $\bb B'\to\bb B$ is meant a right exact functor $G:\bb
B'\to\bb E=\bb E(\bb B,\bb B,\bb B)$ such that $F'=s\circ G$,
$F=t\circ G$, and $F''=q\circ G$.

\begin{prop}[Equivalent formulations of the additivity theorem]\label{eq}

Each of the following conditions implies the three others.

\begin{enumerate}
\item The following projection
   $$i.S.\bb E(\bb B',\bb B,\bb B'')\to i.S.\bb B'\times i.S.\bb B'',\ \ \ E\longmapsto (E_{(0,0)},E_{(1,1)})$$
is a homotopy equivalence.

\item The following projection
   $$i.S.\bb E\to i.S.\bb B\times i.S.\bb B,\ \ \ E\longmapsto (E_{(0,0)},E_{(1,1)})$$
is a homotopy equivalence.

\item The following two maps are homotopic,
   $$i.S.\bb E\to i.S.\bb B,\ \ \ E\longmapsto E_{(0,1)},\textrm{ respectively }E\longmapsto E_{(0,0)}\coprod E_{(1,1)}.$$

\item If $F'\to F\to F''$ is a right exact sequence of right exact functors $\bb B'\to\bb B$ then there exists
a homotopy
   $$|i.S.F|\iso |i.S.F'|\vee|i.S.F''|.$$
\end{enumerate}
\end{prop}

\begin{proof}
(2) is a special case of (1), (3) is a special case of (4), and (4) follows from (3) by naturality.

So it will suffice to show
the implications $(2)\Longrightarrow(3)$ and $(4)\Longrightarrow(1)$.

$(2)\Longrightarrow(3)$. The desired homotopy $|i.S.t|\iso|i.S.(s\vee q)|$ is valid upon the
restriction along the map
   $$|i.S.\bb B|\times|i.S.\bb B|\lra{}|i.S.\bb E|,\ \ \ (A,B)\longmapsto\alpha A\coprod\beta B,$$
so it will suffice to know that this map is a homotopy equivalence. But this map is a section
to the map in~(2) and therefore is a homotopy equivalence if that is one.

$(4)\Longrightarrow(1)$. First consider the maps
$l:\Delta^1\to\Box$, $\kappa\longmapsto(\kappa,0)$, and
$q:\Box\to\Delta^1$, $(\epsilon,\eta)\longmapsto\epsilon$. Denote
by $\bb E'_?=\{q^*l^*E\mid E\in\bb E(\bb B',\bb B,\bb B'')_?\}$.
Given $E\in\bb E(\bb B',\bb B,\bb B'')_0$ the object $q^*l^*E$ is
depicted in $\bb B_0$ as
   $$\begin{CD}
      E_{(0,0)}@>1>>E_{(0,0)}\\
      @VVV@VVV\\
      O@>1>> O
     \end{CD}$$
where $O=E_{(1,0)}$ is a zero object.

Also, let $i:\Delta^1\to\Box$, $\kappa\longmapsto(1,\kappa)$, and
$j:\Box\to\Delta^1$, $(\epsilon,\eta)\longmapsto\eta$. Denote by
$\bb E''_?=\{j^*i^*E\mid E\in\bb E(\bb B',\bb B,\bb B'')_?\}$.
Given $E\in\bb E(\bb B',\bb B,\bb B'')_0$ the object $j^*i^*E$ is
depicted in $\bb B_0$ as
   $$\begin{CD}
      O@>>>E_{(1,1)}\\
      @V1VV@VV1V\\
      O@>>>E_{(1,1)}
     \end{CD}$$

We shall construct a right exact morphism
   $$E\in\bb E(\bb B',\bb B,\bb B'')_?\longmapsto E^2\in\bb E(\bb E',\bb E,\bb E'')_?$$
such that $E^2$ is depicted in $\bb E_0$ as follows.
   $$\begin{CD}
      E'@>>>E\\
      @VVV@VVV\\
      O@>>>E''
     \end{CD}$$
If we neglect the object at $(1,0)\in\Box$ the above diagram is
depicted in $\bb B_0$ as
   $$\begin{CD}
      E_{(0,0)}@>1>>E_{(0,0)}@>>>O\\
      @V1VV@VVV@VVV\\
      E_{(0,0)}@>>>E_{(0,1)}@>>>E_{(1,1)}\\
      @VVV@VVV@VV1V\\
      O@>>>E_{(1,1)}@>>>E_{(1,1)}.\\
     \end{CD}$$
Then it will immediately follow from assumption that
   $$i.S.\bb E(\bb B',\bb B,\bb B'')\to i.S.\bb E'\times i.S.\bb E'',\ \ \ E\longmapsto(q^*l^*E,j^*i^*E)$$
is a homotopy equivalence with section given by
$(E',E'')\longmapsto E'\coprod E''$. We construct the square $E^2$
in the following way. Consider two maps
$\phi,\psi:\Box\to\Delta^1$ defined by
   $$(0,0),(0,1),(1,0)\bl\phi\longmapsto 0,\ \ \ (1,1)\bl\phi\longmapsto 1$$
and
   $$(0,0)\bl\psi\longmapsto 0,\ \ \ (0,1),(1,0),(1,1)\bl\psi\longmapsto 1.$$
Given an object $A\in\bb B_{\Delta^1}$ the functors $\phi^*$ and $\psi^*$ take $A$
to the squares which are depicted in $\bb B_0$ as
   $$\begin{CD}
      A_0@>1>>A_0\\
      @V1VV@VVV\\
      A_0@>>>A_1
     \end{CD}$$
and
   $$\begin{CD}
      A_0@>>>A_1\\
      @VVV@VV1V\\
      A_1@>>1>A_1
     \end{CD}$$
respectively.

An object $E\in\bb E(\bb B',\bb B,\bb B'')_0$ can be regarded as
an object in $\bb B(\Delta^1)_{\Delta^1}$ and it is evaluated as
$E_{(0,0)}\to O$ at $0$ and as $E_{(0,1)}\to E_{(1,1)}$ at 1. We
embedd the $E$ into the cocartesian square
$E^1=(1_{\Delta^1}\times\phi)^*E$ in $\bb B(\Delta^1)_{\Box}$.
After depicting $E^1$ in an appropriate way, one obtains the
following diagram in $\bb B_0$.
   $$\xymatrix@!0{
     &E_{(0,0)}\ar[rr]\ar'[d][dd] && O\ar[dd]\\
     E_{(0,0)}\ar[ur]\ar[rr]\ar[dd] && O\ar[ur]\ar[dd]\\
     &E_{(0,1)}\ar'[r][rr] && E_{(1,1)}\\
     E_{(0,0)}\ar[rr]\ar[ur] && O\ar[ur]
     }$$
The object $E^1$ can be regarded as an object in $\bb
B(\Box)_{\Delta^1}$ and it is evaluated as the left square of the
depicted cube at $0$ and as the right square at 1. We embedd the
$E^1$ into the cocartesian square $E^2=(\psi\times 1_{\Box})^*E^1$
in $\bb B(\Box)_{\Box}$. The construction of $E^2$ is completed.
The morphism $E\in\bb E(\bb B',\bb B,\bb B'')_?\longmapsto
E^2\in\bb E(\bb E',\bb E,\bb E'')_?$ is induced by $(\psi\times
1_{\Box})^*(1_{\Delta^1}\times\phi)^*$.

It remains to show that the maps $f:E\in\bb E'_?\longmapsto
E_{(0,0)}\in\bb B_?$ and $g:E\in\bb E''_?\longmapsto
E_{(1,1)}\in\bb B_?$ induce a homotopy equivalence
   $$i.S.\bb E'\times i.S.\bb E''\xrightarrow{(f,g)}i.S.\bb B'\times i.S.\bb B''.$$
For the map $p:i.S.\bb E(\bb B',\bb B,\bb B'')\to i.S.\bb B'\times
i.S.\bb B''$, $E\longmapsto(E_{(0,0)},E_{(1,1)})$, equals to the
composition
   $$i.S.\bb E(\bb B',\bb B,\bb B'')\xrightarrow{(q^*l^*,j^*i^*)}
     i.S.\bb E'\times i.S.\bb E''\xrightarrow{(f,g)}i.S.\bb B'\times i.S.\bb B''$$
and the left arrow is a homotopy equivalence by above.

Let $\bar{\bb B}_?=\{X\in\bb B_{\Delta^1\times ?}\mid X_1\iso
0\}$. Then $\bar{\bb B}$ and $\bb E'$ are isomorphic because
$l^*q^*=1_{\bar{\bb B}}$ and $q^*l^*q^*l^*|_{\bb E'}=1_{\bb E'}$.
In a similar way, let $\bar{\bar{\bb B}}_?=\{X\in\bb
B_{\Delta^1\times ?}\mid X_0\iso 0\}$. Then $\bar{\bar{\bb B}}$
and $\bb E''$ are isomorphic because $i^*j^*=1_{\bar{\bar{\bb
B}}}$ and $j^*i^*j^*i^*|_{\bb E''}=1_{\bb E''}$.

Finally, the proof of Proposition~\ref{sss} shows that the
morphism $\bar{\bb B}\to\bb B'$ induced by the map $0\longmapsto
0\in\Delta^1$ is an equivalence as well as the morphism
$\bar{\bar{\bb B}}\to\bb B''$ induced by the map $0\longmapsto
1\in\Delta^1$. We are done.
\end{proof}

Given a simplicial object $X$ let $PX\to X$ be the projection
induced by the face map $\partial_0:X_{n+1}\to X_n$. If we
consider $X_1$ as a trivial simplicial object, there is an
inclusion $X_1\to PX$ resulting a sequence $X_1\to PX\to X$.

In particular we obtain a sequence $i.S_1\bb B\to P(i.S.\bb B)\to
i.S.\bb B$ which in view of the equivalence of $i.S_1\bb B$ with
$i.\bb B_0$ may be rewritten as
   $$i.\bb B_0\lra{G} P(i.S.\bb B)\lra{\partial_0}i.S.\bb B.$$

We show explicitly what the map $G$ is. Let $\ell^*:S_1\bb B\to\bb
B_0$ be the equivalence stated in Proposition~\ref{sss}. A
quasi-inverse to $\ell^*$ is constructed as follows. Consider the
open immersion $e:0\longmapsto 0\in\Delta^1$. If $\bb B$ is a left
system of diagram categories and $k:\Delta^1\to 0$ is the map
$0;1\longmapsto 0;\star$, then $(k^*B)_0=B$ and $(k^*B)_1=0$ for
any $B\in\bb B_0$. In turn, if $\bb B$ is a left pointed
d\'erivateur, then $(e_*B)_0=B$ and $(e_*B)_1=0$.

Next, let $p:\Delta^1\to\Ar\Delta^1$ be the closed immersion
$i\longmapsto(i,1)$, $r:\Ar\Delta^{1\star}\to\Delta^{1\star}$ be
the map $(0,0)\longmapsto\star,(0;1,1)\longmapsto 0;1$, and
$B\in\bb B_0$. Then $(r^*k^*B)_{(0,0)}=(r^*k^*B)_{(1,1)}=0$ and
$(r^*k^*B)_{(0,1)}=B$ if $\bb B$ is a left system of diagram
categories and denote $g=r^*k^*$. If $\bb B$ is a left pointed
d\'erivateur, it follows that
$(p_!e_*B)_{(0,0)}=(p_!e_*B)_{(1,1)}=0$ and $(p_!e_*B)_{(0,1)}=B$.
In this case $g=p_!e_*$.

Consider further the map $v:\Delta^{n+1}\to\Delta^1$,
$0\longmapsto 0$ and $i\longmapsto 1$ for $i\geq 1$. We put
$G=v^*g:\bb B_0\to S_{n+1}\bb B$. Then the ``values" of $GB$ at
each $(i,j)\in\Ar\Delta^{n+1}$ are: $GB_{(i,j)}=0$ if
$(i,j)=(0,0)$ and if $i\geq 1$ and $GB_{(0,j)}=B$ for any $j\geq
1$. Considering $\bb B_0$ as a trivial simplicial category, there
results the map $G:\bb B_0\to PS.\bb B$ as well as the map
$G:|i.\bb B_0|\to|P(i.S.\bb B)|$.

The composite map $|i.\bb B_0|\lra{G}|P(i.S.\bb B)|\to|i.S.\bb B|$
is constant, and $|P(i.S.\bb B)|$ is contractible (for it is
homotopy equivalent to the contractible space $|i.S_0\bb B|$), so
we obtain a map, well defined up to homotopy,
   $$|i.\bb B_0|\to\Omega|i.S.\bb B|.$$

We make a couple of useful observations due to
Waldhausen~\cite[p.~332]{W}.

\begin{observation}
The following two composite maps are homotopic,
   $$\xymatrix{|i.\bb E_0|\ar@<2.5pt>[r]^t
                        \ar@<-2.5pt>[r]_{s\vee q}&|i.\bb B_0|\ar[r]&\Omega|i.S.\bb B|.}$$
\end{observation}

\begin{proof}
This results from an inspection of $|i.S.\bb B|_{(2)}$, the 2-skeleton of $|i.S.\bb B|$
in the $S.$-direction. We can identify $i\bb B_0$ to $iS_1\bb B$ and $i\bb E_0$
to $iS_2\bb B$.

The face maps from $i.S_2\bb B$ to $i.S_1\bb B$ then correspond to
the three maps $s,t,q$, respectively, and each of which can be
seen from the diagram
   $$\xymatrix{0\ar[rr]^{A_{0,2}} \ar[dr]_{A_{0,1}} && 2 &\\
               &1\ar[ur]_{A_{1,2}} &}$$
Let us consider the canonical map $|i.S_2\bb B|\times|\Delta^2|\to|i.S.\bb B|_{(2)}$.
Regarding the 2-simplex $|\Delta^2|$ as a homotopy from the edge $(0,2)$ to the edge
path $(0,1)(1,2)$ we obtain a homotopy from the composite map $jt$,
   $$|i.\bb E_0|\lra{t}|i.\bb B_0|\lra{j}\Omega|i.S.\bb B|_{(2)}$$
to the loop product of two composite maps $js$ and $jq$. But in $\Omega|i.S.\bb B|$
the loop product is homotopic to the composition law, by a well known fact about
loop spaces of $H$-spaces, whence the observation as stated.
\end{proof}

The same consideration shows, more generally,

\begin{observation}
For every $n\geq 0$ the two composite maps
   $$\xymatrix{|i.S.^n\bb E|\ar@<2.5pt>[r]^t
                        \ar@<-2.5pt>[r]_{s\vee q}&|i.S.^n\bb B|\ar[r]&\Omega|i.S.^{n+1}\bb B|.}$$
are homotopic.
\end{observation}

\begin{thm}\label{does}
The additivity theorem (Proposition~\ref{eq}) is valid if the
definition of $K$-theory as $\Omega|i.S.\bb B|$ is substituted
with $\Omega^\infty|i.S.^\infty\bb B|= \lim_n\Omega^n|i.S.^n\bb
B|$.
\end{thm}

\begin{proof}
By the preceding observation the two composite maps
   $$\xymatrix{\Omega^\infty|i.S.^\infty\bb E|\ar@<2.5pt>[r]^t\ar@<-2.5pt>[r]_{s\vee q}
     & \Omega^\infty|i.S.^\infty\bb B|\ar[r]
     &\Omega^\infty|i.S.^\infty\bb B|}$$
are homotopic. Since the map on the right is an isomorphism, this is one of the
equivalent conditions of the additivity theorem (Proposition~\ref{eq}).
\end{proof}

\begin{rem}{\rm
As a consequence of the theorem we could add yet another
reformulation of the additivity theorem to the list of
Proposition~\ref{eq} (see also Theorem~\ref{jjj}). Namely the
additivity theorem as stated there implies that the maps
$|i.S.^n\bb B|\to\Omega|i.S.^{n+1}\bb B|$ are homotopy
equivalences for $n\geq 1$. Conversely if these maps are homotopy
equivalences then so is $\Omega|i.S.\bb
B|\to\Omega^\infty|i.S.^\infty\bb B|$, and thus the additivity
theorem is provided by the theorem. }\end{rem}

Let $F:\bb A\to\bb B$ be a right exact functor between two left
systems of diagram categories or between two left pointed
d\'erivateurs respectively. We denote by $\bb S.(F:\bb A\to\bb B)$
the fibred product of the diagram
   $$\bb S.\bb A\lra F\bb S.\bb B\bl{\partial_0}\longleftarrow P\bb S.\bb B$$
with $\partial_0=d_0^*$ the map induced by
$d_0:\Delta^n\to\Delta^{n+1}$. By Proposition~\ref{pp} $\bb
S.(F:\bb A\to\bb B)$ is a simplicial left system of diagram
categories or a simplicial left pointed d\'erivateur respectively.
Thus for every $n$ one has a commutative diagram
   $$\begin{CD}
      \bb S_n(F:\bb A\to\bb B)@>{F'}>>(P\bb S.\bb B)_n=\bb S_{n+1}\bb B\\
      @VpVV@VV{\partial_0}V\\
      \bb S_n\bb A@>F>>\bb S_n\bb B.
     \end{CD}$$
By construction we can identify an object of $\bb S_n(F:\bb
A\to\bb B)_?$ to a triple $(A,c,B)$ with objects in $\bb S_n\bb
A_?$ and $\bb S_{n+1}\bb B_?$, respectively, together with an
isomorphism $FA\bl c\iso\partial_0B$. We note that all the maps in
the defined diagram are right exact.

Let $G:\bb B\to\bb S_{n+1}\bb B$ be the morphism constructed
above. We have $\partial_0GB=0$. Then $G$ factors as $F'\circ G'$
with $G':\bb B\to\bb S_n(F:\bb A\to\bb B)$,
$B\bl{G'}\longmapsto(0,1,GB)$.

Regarding $\bb B$ as a simplicial object in a trivial way, we obtain a sequence
   \begin{equation}\label{qw}
    \bb B\lra{G'}\bb S.(F:\bb A\to\bb B)\lra{p}\bb S.\bb A
   \end{equation}
in which the composed map is trivial. There results a sequence
   $$i.S.\bb B\lra{}i.S.S.(\bb A\to\bb B)\lra{}i.S.S.\bb A,$$
induced by~\eqref{qw}.

Similarly, there is a sequence
   $$i.S.\bb B\to P(i.S.S.\bb B)\to i.S.S.\bb B$$
where the ``$P$" refers to the first $S.$-direction, say.

\begin{thm}\label{jjj}
The following statements are equivalent:

$(1)$ the additivity theorem (Proposition~\ref{eq}) is valid;

$(2)$ the sequence
   $$i.S.\bb B\lra{}i.S.S.(\bb A\to\bb B)\lra{}i.S.S.\bb A,$$
is a fibration up to homotopy;

$(3)$ the sequence
   $$i.S.\bb B\to P(i.S.S.\bb B)\to i.S.S.\bb B$$
is a fibration up to homotopy;

$(4)$ the map $|i.S.^n\bb B|\to\Omega|i.S.^{n+1}\bb B|$ is a homotopy equivalence
for any $n\geq 1$.

If the equivalent conditions $(1)-(4)$ hold, then the spectrum
   $$n\longmapsto i.S.^n\bb B$$
with structural maps being defined as the map $|i.\bb
B_0|\to\Omega|i.S.\bb B|$ above is a $\Omega$-spectrum beyond the
first term. The spectrum is connective (the $n$th term is
$(n-1)$-connected). As a consequence, the $K$-theory for $\bb B$
can then equivalently be defined as the space
   $$\Omega^\infty|i.S.^\infty\bb B|=\lim_n\Omega^n|i.S.^n\bb B|.$$
\end{thm}

\begin{proof}
(3) is a consequence of (2). Since the space $|P(i.S.S.\bb B)|$ is
contractible, the condition (3) implies that the map $|i.S.\bb
B|\to\Omega|i.S.S.\bb B|$ is a homotopy equivalence and more
generally therefore also the map $|i.S.^n\bb
B|\to\Omega|i.S.^{n+1}\bb B|$ for any $n\geq 1$. So (4) is a
consequence of (3). By the second observation following
Proposition~\ref{eq} the two composite maps
   $$\xymatrix{|i.S.\bb E|\ar@<2.5pt>[r]^t
                        \ar@<-2.5pt>[r]_{s\vee q}&|i.S.\bb B|\ar[r]&\Omega|i.S.S.\bb B|}$$
are homotopic. If the map on the right is a homotopy equivalence then $t$ is homotopic to
$s\vee q$. Thus (4) implies~(1). It remains therefore to prove $(1)\Longrightarrow(2)$.

By Lemma~\ref{wald} it is enough to show that for every $n$ the
sequence $i.S.\bb B\lra{}i.S.S_n(\bb A\to\bb B)\lra{}i.S.S_n\bb A$
is a fibration (since the base term $i.S.S_n\bb A$ is connected
for every $n$). Using the additivity theorem we will show that the
sequence is the same, up to homotopy, as the trivial fibration
sequence associated to the product $i.S.\bb B\times i.S.S_n\bb A$.

Consider the maps $u:\Delta^1\to\Delta^{n+1}$, $0;1\longmapsto
0;1$, and $v:\Delta^{n+1}\to\Delta^1$, $0\longmapsto 0$,
$i\longmapsto 1$ for $i\geq 1$. Then $u$ is left adjoint to $v$.
To simplify the notation the corresponding maps
$\Ar\Delta^1\to\Ar\Delta^{n+1}$ and
$\Ar\Delta^{n+1}\to\Ar\Delta^1$ induced by $u$ and $v$ denote by
the same letters. Let $\bar{\bb B}=\{v^*u^*B\mid B\in\bb
S_{n+1}\bb B\}$. It follows that $v^*u^*B_{(0,0)}=B_{(0,0)}$,
$v^*u^*B_{(0,i)}=B_{(0,1)}$ for any $1\le i\le n+1$, and
$v^*u^*B_{(i,j)}=B_{(1,1)}$ for any $i\geq 1$.

Denote by $\bar{\bar{\bb B}}=\{\sigma_0\partial_0B\mid B\in\bb
S_{n+1}\bb B\}$, where $\sigma_0:\bb S_n\bb B\to\bb S_{n+1}\bb B$
is the map induced by $s_0:\Ar\Delta^{n+1}\to\Ar\Delta^n$. Note
that $\sigma_0$ is right adjoint to $\partial_0$.

Let $m:\Delta^1\times\Ar\Delta^{n+1}\to\Ar\Delta^{n+1}$ be the map
taking $(0,(i,j))$ to $(uv(i),uv(j))$ and $(1,(i,j))$ to $(i,j)$.
Then $m^*$ takes an object $B\in\bb S_{n+1}\bb B_?$ to that in
$\bb S_{n+1}\bb B_{\Delta^1\times ?}$, which is depicted in $\bb
S_{n+1}\bb B_?$ as the adjunction morphism $v^*u^*B\to B$.

Next, let $l:\Delta^1\times\Ar\Delta^{n+1}\to\Ar\Delta^{n+1}$ be
the map taking $(0,(i,j))$ to $(i,j)$ and $(1,(i,j))$ to
$(d_0s_0(i),d_0s_0(j))$. Then $l^*$ takes an object $B\in\bb
S_{n+1}\bb B_?$ to the object in $\bb S_{n+1}\bb B_{\Delta^1\times
?}$, which is evaluated as $B$ at 0, and as $\sigma_0\partial_0B$
at 1. This object is depicted in $\bb S_{n+1}\bb B_?$ as the
adjunction morphism $\beta:B\to\sigma_0\partial_0B$.

Restriction of the right exact morphism $(1_{\Delta^1}\times
m)^*l^*:\bb B(\Ar\Delta^{n+1})\to \bb B(\Box\times\Delta^{n+1})$
to $\bb S_{n+1}\bb B$ takes an object $B\in\bb S_{n+1}\bb B_?$ to
one in $\bb E(\Ar\Delta^{n+1})_?\subset\bb
B(\Box\times\Ar\Delta^{n+1})_?$. Thus we result in a right exact
functor
   $$T:\bb S_{n+1}\bb B\lra{l^*}\bb S_{n+1}\bb B(\Delta^1)\xrightarrow{(1_{\Delta^1}\times m)^*}\bb E(\Ar\Delta^{n+1})$$
such that $t\circ T$ is the identity morphism on $\bb S_{n+1}\bb
B$ and for every $B\in\bb S_{n+1}\bb B$ we have $s\circ
T(B)\in\bar{\bb B}$ and $q\circ T(B)\in\bar{\bar{\bb B}}$. So $T$
takes its values in $\bb E(\bar{\bb B},\bb S_{n+1}\bb
B,\bar{\bar{\bb B}})$ and it is actually a functor
   $$T:\bb S_{n+1}\bb B\to\bb E(\bar{\bb B},\bb S_{n+1}\bb B,\bar{\bar{\bb B}}).$$

To illustrate the reader the above procedure, let us think of $\bb
S_{n+1}\bb B$ for a short while as ``strings" $\bb B(\Delta^n)$
via the equivalence $\ell^*:\bb S_{n+1}\bb B\to\bb B(\Delta^n)$
stated in Proposition~\ref{sss}. Let us consider the function $\wt
m:\Delta^1\times\Delta^n\to\Delta^n$ defined by
   $$(0,i)\longmapsto 0,\ \ \ (1,i)\longmapsto i.$$
Then the induced right exact morphism $\wt m^*:\bb
B(\Delta^n)\to\bb B(\Delta^1\times\Delta^n)$ takes an object
$B\in\bb B(\Delta^n)_?$ to the object in $\bb
B(\Delta^1\times\Delta^n)_?$ which is depicted in $\bb B_?$ as
   $$\begin{CD}
      B_0@>1>>B_0@>1>>\cdots@>1>>B_0\\
      @V1VV@VV{b_1}V@.@VV{b_{n}\cdots b_1}V\\
      B_0@>{b_1}>>B_1@>{b_2}>>\cdots@>{b_{n}}>>B_n
     \end{CD}$$
Let $k:\Delta^1\times\Delta^n\to\Ar\Delta^{n+1}$ be the map
$(i,j)\longmapsto(i,j+1)$,
$\alpha_j:\Box\to\Delta^1\times\Delta^n$ the map taking $(0;1,0)$
to $(0;1,j)$ and $(0;1,1)$ to $(0;1,j+1)$. Denote by $\bb
S_{n+1}'\bb B$ the left subsystem of diagram categories or the
left pointed subd\'erivateur of $\bb B(\Delta^1\times\Delta^n)$
respectively consisting of the objects $B$ such that all the
squares $\alpha_j^*B$, $j\le n$, are cocartesian and $B_{(1,0)}=O$
is a zero object. Then the restriction morphisms $k^*:\bb
S_{n+1}\bb B\to\bb S_{n+1}'\bb B$ and $w^*:\bb S_{n+1}'\bb B\to\bb
B(\Delta^n)$ with $w:\Delta^n\to\Delta^1\times\Delta^n$,
$j\longmapsto(0,j)$ are equivalences as one easily shows.

Restriction of the morphism $(1_{\Delta^1}\times \wt m)^*:\bb
B(\Delta^1\times\Delta^n)\to \bb B(\Box\times\Delta^n)$ to $\bb
S_n'\bb B$ takes an object $B\in\bb S_n'\bb B_?$ to one in $\bb
B(\Box\times\Delta^n)_?$ which is depicted in $\bb B_?$ as
   $$\xymatrix@!0{
     &B_{(0,0)}\ar[rr]\ar'[d][dd] && B_{(0,1)}\ar'[d][dd] \ar[rr]\ar'[d][dd] && B_{(0,2)}\ar'[d][dd] \ar[rr]\ar'[d][dd]
     && \cdots \ar[rr] && B_{(0,n)}\ar[dd]\\
     B_{(0,0)}\ar[ur]\ar[rr]\ar[dd] && B_{(0,0)}\ar[ur]\ar[dd] \ar[ur]\ar[rr]\ar[dd]
     && B_{(0,0)}\ar[ur]\ar[dd] \ar[ur]\ar[rr]\ar[dd]&& \cdots\ar[rr] &&B_{(0,0)}\ar[ur]\ar[dd]\\
     &O\ar'[r][rr] && B_{(1,1)} \ar'[r][rr] && B_{(1,2)} \ar[rr] && \cdots \ar'[r][rr] && B_{(1,n)}\\
     O\ar[rr]\ar[ur] && O\ar[ur]\ar[rr] && O\ar[ur] \ar[rr]&& \cdots \ar[rr] && O\ar[ur]
     }$$
The back wall of the diagram is the element $B$ of $\bb S_n'\bb
B_?$ drawn in $\bb B_?$. We obtain a right exact morphism $\wt
T=(1_{\Delta^1}\times \wt m)^*\circ w^*{}^{-1}:\bb
B(\Delta^n)\to\bb E(\Delta^n)$. For each transverse $x$th square
$B_{x,\Box}$, $x\le n$, is cocartesian.

Now, return to the morphism $T$ and lift it to $\bb S_n(\bb
A\to\bb B)$. To be more precise, we send an object $(A,c,B)\in\bb
S_n(\bb A\to\bb B)_?$ to
$(\partial_0T\sigma_0A,\partial_0T\sigma_0(c),TB)\in\bb E(\bb
S_n(\bb A\to\bb B))_?$. We use the commutative diagram
   $$\xymatrix{
      \bb A(\Ar\Delta^n)\ar[r]^{\sigma_0}\ar[d]_F&\bb
      A(\Ar\Delta^{n+1})\ar[r]^{T}\ar[d]_F&\bb A(\Box\times\Ar\Delta^{n+1})\ar[r]^{\partial_0}\ar[d]^F
      &\bb A(\Box\times\Ar\Delta^n)\ar[d]^{F}\\
      \bb B(\Ar\Delta^n)\ar[r]^{\sigma_0}&\bb
      B(\Ar\Delta^{n+1})\ar[r]^{T}&\bb B(\Box\times\Ar\Delta^{n+1})\ar[r]^{\partial_0}
      &\bb B(\Box\times\Ar\Delta^n)
     }$$
to show that $F\partial_0T\sigma_0A=\partial_0T\sigma_0FA$. The
relation $\partial_0TB=\partial_0T\sigma_0\partial_0B$ is
straightforward.

Let $\bb
B'=\{(\partial_0v^*u^*\sigma_0A,\partial_0v^*u^*\sigma_0(c),v^*u^*B)\mid(A,c,B)\in\bb
S_n(\bb A\to\bb B)\}$ and $\bb
B''=\{(A,c,\sigma_0\partial_0B)\mid(A,c,B)\in\bb S_n(\bb A\to\bb
B)\}$.

There results a right exact functor
   $$T':\bb S_n(\bb A\to\bb B)\to\bb E(\bb B',\bb S_n(\bb A\to\bb B),\bb B'')$$
with $s\circ T'$ sending $(A,c,B)$ to
$(\partial_0v^*u^*\sigma_0A,\partial_0v^*u^*\sigma_0(c),v^*u^*B)$,
$t\circ T'$ being the identity, and $q\circ T'$ sending $(A,c,B)$
to $(A,c,\sigma_0\partial_0B)$. Thus we get an exact sequence
$s\circ T'\to 1\to q\circ T'$. It follows from our assumption that
the map
   $$(s\circ T',q\circ T'):S.S_n(\bb A\to\bb B)\to S.\bb B'\times S.\bb B''$$
is a homotopy equivalence with a homotopy inverse induced by
coproduct.

Clearly, the morphism $\bb B'\to\bb B$ taking
$(\partial_0v^*u^*\sigma_0A,\partial_0v^*u^*\sigma_0(c),v^*u^*B)$
to $B_{(0,1)}$ is an equivalence. Its quasi-inverse is given by
$G'$.

Let us show that the morphism $\delta:\bb S_n\bb A\to\bb B''$,
$A\longmapsto(A,1,\sigma_0FA)$, is a quasi-inverse to the
restriction of $p$ to $\bb B''$. Obviously $\delta$ is faithful.
Given an object $(A,c,B)\in\bb B''$ the map
$(1,\sigma_0(c)):(A,1,\sigma_0FA)\to(A,c,B)$ is an isomorphism. It
also follows that every map $(a,b):\delta A\to\delta A'$ in $\bb
B''$ equals to $(a,\sigma_0Fa)$, and hence $\delta$ is also full.
We see that $\delta$ is an equivalence.

It follows that the map
   $$i.S.\bb B'\times i.S.\bb B''\to i.S.\bb B\times i.S.S_n\bb A$$
is a homotopy equivalence, hence is so the composite
   $$i.S.S_n(\bb A\to\bb B)\to i.S.\bb B'\times i.S.\bb B''\to i.S.\bb B\times i.S.S_n\bb A.$$
This homotopy equivalence fits into the following commutative diagram
   $$\begin{CD}
      i.S.\bb B@>>>i.S.S_n(\bb A\to\bb B)@>>>i.S.S_n\bb A\\
      @V1VV@VVV@VV1V\\
      i.S.\bb B@>>>i.S.\bb B\times i.S.S_n\bb A@>>>i.S.S_n\bb A
     \end{CD}$$
Being homotopy equivalent to the trivial fibration (the lower row
of the diagram), we conclude that the upper sequence is a
fibration, as was to be shown.
\end{proof}

\begin{rem}\label{jjjj}{\rm
Let $\Im$ be either a class of left systems of diagram categories
or left pointed d\'erivateurs satisfying the following two
conditions:

(1) $\bb B\in\Im$ implies $\bb S_n\bb B\in\Im$ for any $n$;

(2) the map $i.S.\bb E\xrightarrow{(s,q)}i.S.\bb B\times i.S.\bb
B$ is a homotopy equivalence for any $\bb B\in\Im$.

The proof of Theorem~\ref{jjj} then shows that the spectrum
   $$n\longmapsto i.S.^n\bb B$$
is a $\Omega$-spectrum beyond the first term, and so the
$K$-theory for every $\bb B\in\Im$ can then equivalently be
defined as the space
   $$\Omega^\infty|i.S.^\infty\bb B|=\lim_n\Omega^n|i.S.^n\bb B|.$$
}\end{rem}

A left pointed d\'erivateur $\bb D$ of the domain $\ord$ is said
to be {\it complicial\/} if there is a right exact equivalence
$F:\bb D\cc C\lra{}\bb D$ for some complicial biWaldhausen
category $\cc C$ in the sense of Thomason and which is closed
under the formation of canonical homotopy pushouts and canonical
homotopy pullbacks. In this case we say that $\bb D$ is {\it
represented\/} by $\cc C$. That equivalence induces a homotopy
equivalence of bisimplicial sets $F:i.S.\bb D\cc C\lra{}i.S.\bb
D$.

\begin{thm}[\cite{Gar}]\label{complic}
The class of complicial d\'erivateurs satisfies the conditions of
the remark above.
\end{thm}

\begin{prop}\label{111}
Under the hypotheses of Theorem~\ref{jjj} suppose we are given a
sequence $\bb A\to\bb B\to\bb C$ of right exact morphisms between
left systems of diagram categories or left pointed d\'erivateurs
respectively. Then the square
   $$\begin{CD}
      i.S.\bb B@>>>i.S.S.(\bb A\to\bb B)\\
      @VVV@VVV\\
      i.S.\bb C@>>>i.S.S.(\bb A\to\bb C)
     \end{CD}$$
is homotopy cartesian.
\end{prop}

\begin{proof}
There is a commutative diagram
   $$\begin{CD}
      i.S.\bb B@>>>i.S.S.(\bb A\to\bb B)@>>>i.S.S.\bb A\\
      @VVV@VVV@VV{\id}V\\
      i.S.\bb C@>>>i.S.S.(\bb A\to\bb C)@>>>i.S.S.\bb A
     \end{CD}$$
in which the rows are fibrations up to homotopy by
Theorem~\ref{jjj}. Therefore the square on the left is homotopy
cartesian.
\end{proof}

\begin{cor}\label{333}
Under the hypotheses of Theorem~\ref{jjj} the following two assertions are valid.

$(1)$ To a right exact morphism there is associated a fibration
   $$i.S.\bb B\to i.S.\bb C\to i.S.S.(\bb B\to\bb C).$$

$(2)$ If $\bb C$ is a retract of $\bb B$ (by right exact functors) there is a splitting
   $$i.S.\bb B\iso i.S.\bb C\times i.S.S.(\bb C\to\bb B).$$
\end{cor}

\begin{proof}
(1). If $\bb A=\bb B$ the space $|i.S.S.(\bb A=\bb A)|$ is contractible whence the first assertion.

(2). This is the case of Proposition~\ref{111} where the composed map $\bb A\to\bb B\to\bb C$
is an identity map since $i.S.S.(\bb A\to\bb C)$ is contractible in that case.
\end{proof}

\section{Concluding remarks}

Given an exact category $\cc E$, one would like to compare
Quillen's $K$-theory $K(\cc E)$ of $\cc E$ with the $K$-theory of
the associated bid\'erivateur $\bb D^b(\cc E)$.

Let $wC^b(\cc E)$ denote the Waldhausen category of
quasi-isomorphisms in $C^b(\cc E)$ with cofibrations componentwise
admissible monomorphisms. We have a natural functor for every
$I\in\dirf$,
   $$\Ho:C^b(\cc E^I)\to D^b(\cc E^I).$$
The image under the functor $\Ho$ of any cocartesian square of
$C^b(\cc E)^\Box=C^b(\cc E^\Box)$
   $$\begin{CD}
      *@>>>*\\
      @VVV@VVV\\
      *@>>>*
     \end{CD}$$
in which the horizontal arrows are cofibrations is a cocartesian
square in $\bb D^b(\cc E)_\Box$ (this is dual to~\cite[3.14]{C1}).
Therefore $\Ho$ induces a map of bisimplicial objects
   $$\nu:w.S.C^b(\cc E)\to i.S.\bb D^b(\cc E).$$
Consider the map
   $$K(\tau):K(\cc E)\to K(wC^b(\cc E))$$
which is induced by the map $\tau$ taking an object of $\cc E$ to
the complex concentrated in the zeroth degree ($K(wC^b(\cc E))$
stands for the Waldhausen $K$-theory of $wC^b(\cc E))$.

\begin{question}[The first Maltsiniotis conjecture~\cite{M}]
The map $K(\rho)=K(\nu\tau):K(\cc E)\to K(\bb D^b(\cc E))$ is a
homotopy equivalence.
\end{question}

The particular map $K_0(\cc E)\to K_0(\bb D^b(\cc E))$ is an
isomorphism for the Grothen\-dieck groups $K_0(\cc E)$ and
$K_0(D^b(\cc E))$ are naturally isomorphic (exercise!) and
$K_0(\bb D^b(\cc E))$ is naturally isomorphic to $K_0(D^b(\cc E))$
by Lemma~\ref{555}.


The first Malsiniotis conjecture is very resistant in general.
However one can obtain some information for a large class of exact
categories including the abelian categories. The following shows
that Quillen's $K$-theory $K(\cc E)$ of an exact category $\cc E$
from this class is a retract of $K(\bb D^b(\cc E))$.

\begin{thm}\label{neeman}
Let $\cc E$ be an extension closed full exact subcategory of an
abelian category $\cc A$ satisfying the conditions of the
Resolution Theorem. That is

$(1)$ if $\es {M'}M{M''}$ is exact in $\cc A$ and $M,M''\in\cc E$,
then $M'\in\cc E$ and

$(2)$ for any object $M\in\cc A$ there is a finite resolution
$0\to P_n\to P_{n-1}\to\cdots\to P_0\to M\to 0$ with $P_i\in\cc
E$.

Then the map
   $$K(\rho):K(\cc E)\to K(\bb D^b(\cc E))$$
is a split inclusion in homotopy. There is a map
   $$p:K(\bb D^b(\cc E))\to K(\cc E)$$
which is left inverse to it. That is $p\circ K(\rho)$ is homotopic
to the identity. In particular, each $K$-group $K_n(\cc E)$ is a
direct summand of $K_n(\bb D^b(\cc E))$.
\end{thm}

We postpone the proof. One should remark that it essentially uses
Neeman's results~\cite{N} on $K$-theory for triangulated
categories.

It is shown in~\cite{To1,TV} that the natural morphism $K(\cc
C)\to K(\bb D\cc C)$ from the Waldhausen $K$-theory to the
$K$-theory of its d\'erivateur can not be an equivalence in
general. For instance this is so for the Waldhausen $K$-theory of
spaces. This does not mean however that the $K$-groups $K_n(\cc
C)$ can not be reconstructed from its d\'erivateur and that this
is a counter example to the comparison problem above stated for
exact categories. The obstruction is really concerned with
functoriallity at the level of spectra.

A {\it good\/} Waldhausen category is a Waldhausen category that
can be embedded in the category of cofibrant objects of a pointed
model category, and whose Waldhausen structrure is induced by the
ambient model structure (see the precise definition in~\cite{TV}).
Though there exist non-good Waldhausen categories
(see~\cite[Example~2.2]{TV}), in practice it turns out that given
a Waldhausen category there is always a good Waldhausen model,
i.e. a good Waldhausen category with the same $K$-theory space up
to homotopy. Any good Waldhausen category is a Waldhausen category
of cofibrant objects, and therefore one can associate to it the
left pointed d\'erivateur $\bb D\cc C$ (Theorem~\ref{cis}). The
following  theorem is also a consequence of a result by Cisinski
and To\"en~\cite[2.16]{To1}.

\begin{thm}[first stated by To\"en~\cite{To}]\label{toen}
Let $\cc C$ and $\cc C'$ be two good Waldhausen categories such
that their associated d\'erivateurs $\bb D\cc C$ and $\bb D\cc C'$
are equivalent. Then the Waldhausen $K$-theory spectra $K(\cc C)$
and $K(\cc C')$ are equivalent as well.
\end{thm}

Under certain extra data, $\bb B$ encodes the structure of a
triangulated category on $\bb B_0$~\cite{C1,F,M1,M}. This
structure is canonically carried over all the categories $\bb
B_I$, $I\in\di$. In this case $\bb B$ is referred to as a {\it
system of triangulated diagram categories\/} or {\it triangulated
d\'erivateur\/} respectively. The following result shows that such
$\bb B$ contains strictly more information than its triangulated
category $\bb B_0$.

\begin{prop}
There exist two non-equivalent triangulated d\'erivateurs $\bb B$
and $\bb B'$, whose associated triangulated categories $\bb B_0$
and $\bb B_0'$ are equivalent.
\end{prop}

\begin{proof}
Let $\cc C=m\cc M(\bb Z/p^2)$ and $\cc C'=m\cc M(\bb
Z/p[\varepsilon]/\varepsilon^2)$ be two stable model categories
considered in~\cite{Sch}. Here $\cc M(\bb Z/p^2)$ and $\cc M(\bb
Z/p[\varepsilon]/\varepsilon^2)$ denote the corresponding
categories of finitely generated modules. Since both $\bb Z/p^2$
and $\bb Z/p[\varepsilon]/\varepsilon^2$ are quasi-Frobenius
rings, it follows that $\cc M(\bb Z/p^2)$ and $\cc M(\bb
Z/p[\varepsilon]/\varepsilon^2)$ are Frobenius categories and $\bb
D\cc C$ and $\bb D\cc C'$ are triangulated d\'erivateurs
by~\cite[4.19]{C1}. It follows from~\cite[1.4]{Sch} that $\bb D\cc
C_0$ and $\bb D\cc C'_0$ are equivalent as triangulated
categories. But the d\'erivateurs $\bb D\cc C$ and $\bb D\cc C'$
can not be equivalent by Theorem~\ref{toen}, because the
Waldhausen $K$-theories $K(\cc C)$ and $K(\cc C')$ are not
equivalent by~\cite[1.7]{Sch}.
\end{proof}

Another problem arising in our context (see
also~\cite[Conjecture~2]{M}) is the localization theorem. Suppose
we are given a family $\cc W=\{\cc W_I\subseteq\Mor\bb B_I\mid
I\in\di\}$ of morphisms compatible with the structure functors
$f^*$ and $f_!$; that is $f^*(\cc W_J)\subseteq\cc W_I$ and
$f_!(\cc W_I)\subseteq\cc W_J$ for every map $f:I\to J$. Let $\bb
B_?[\cc W^{-1}_?]$ denote the category of fractions obtained by
inverting the maps in $\cc W_?$. We also require the following
condition to hold: a morphism is in $\cc W_?$ iff its image in
$\bb B_?[\cc W^{-1}_?]$ is an isomorphism. Let the hyperfunctor
    $$I\bl Q\longmapsto\bb B_I[\cc W^{-1}_I]$$
determine a left system of diagram categories or a left pointed
d\'erivateur respectively. Denote it by $\bb B[\cc W^{-1}]$.
Suppose further that the quotient morphism $Q:\bb B\to\bb B[\cc
W^{-1}]$ is right exact.

In the case when $\bb B$ is a system of triangulated diagram
categories or a triangulated d\'erivateur respectively, then any
thick subcategory $\bb A_0$ of $\bb B_0$ gives rise to a
localization in $\bb B$. Precisely, given $I\in\di$ let $\bb
A_I=\{A\in\bb B_I\mid A_x\in\bb A_0\textrm{ for all } x\in I\}$.
Then $\bb A_I$ is thick in $\bb B_I$ and the functor
   $$I\longmapsto\bb A_I$$
determines a system of triangulated diagram categories or a
triangulated d\'erivateur respectively and the quotient is then
naturally constructed (see~\cite[p.~39]{F}).

\begin{question}[The second Maltsiniotis conjecture~\cite{M}]
Suppose we are given a sequence of morphisms between left systems
of diagram categories or left pointed d\'erivateurs respectively,
   $$\bb A\lra{F}\bb B\lra{Q}\bb B[\cc W^{-1}]$$
where $Q$ is the quotient morphism and $F$ is a right exact
equivalence between $\bb A$ and $Q^{-1}(0)=\{B\in\bb B_?\mid 0\to
B\in\cc W_?\}$. Then the induced sequence of $K$-theory spaces
   $$K(\bb A)\lra{}K(\bb B)\lra{}K(\bb B[\cc W^{-1}])$$
is a fibration up to homotopy.
\end{question}

We have already associated to the morphism $F$ a fibration (see
Corollary~\ref{3333}(1))
   $$i.S.\bb A\to i.S.\bb B\to i.N.S.(\bb A\to\bb B).$$
There is a natural map from $i.N.S.(\bb A\to\bb B)$ to $i.S.\bb
B[\cc W^{-1}]$. Therefore the localization theorem is reduced,
say, to showing that the latter map is a homotopy equivalence.

To conclude, we should also mention another natural construction
associated to a model category $\cc C$, the simplicial
localization $L^H\cc C$, which should carry roughly the same
homotopical information about $\cc C$ as its d\'erivateur $\bb
D\cc C$. Given a good Waldhausen category $\cc C$, To\"en and
Vezzosi~\cite{TV} associate to $L^H\cc C$ a $K$-theory space
$K(L^H\cc C)$ and show that the Waldhausen $K$-theory $K(\cc C)$
is equivalent to $K(L^H\cc C)$. We also recommend the reader to
consult To\"en's thesis~\cite{To1}.

It remains to prove, as promised, Theorem~\ref{neeman}. We start
with preparations.

\begin{defs}{\rm
An additive category $\cc T$ will be called a {\it category with
squares\/} provided

\begin{enumerate}
\item[$\diamond$] $\cc T$ has an automorphism $\Sigma:\cc T\to\cc
T$;

\item[$\diamond$] $\cc T$ comes equipped with a collection of {\it
special squares\/}
\end{enumerate}
   $$\xymatrix{C\ar[r]&D\ar@/^50pt/[dl]^{(1)}\\
               A\ar[u]\ar[r]&B\ar[u]}$$
This means that the square
   $$\xymatrix{C\ar[r]&D\\
               A\ar[u]\ar[r]&B\ar[u]}$$
is commutative in $\cc T$, and there is a map $D\to\Sigma A$
depicted as the curly arrow. The (1) in the label of the arrow is
to remind us that the map is of degree 1, that is a map
$D\to\Sigma A$.

Given two categories with squares, a {\it special functor}
   $$F:\cc S\to\cc T$$
is an additive functor such that there is a natural isomorphism
$\Sigma F\iso F\Sigma$ and $F$ takes special squares in $\cc S$ to
special squares in $\cc T$.

If $\cc T$ is a category with squares, the {\it fold of the
square}
   $$\xymatrix{C\ar[r]^{\delta}&D\ar@/^50pt/[dl]^{\mu}\\
               A\ar[u]^{\beta}\ar[r]^{\alpha}&B\ar[u]^{\gamma}}$$
will be the sequence
   $$A\xrightarrow{(\alpha,-\beta)^t}B\ps C\xrightarrow{(\gamma,\delta)}D\lra{\mu}\Sigma A.$$

}\end{defs}

\begin{exs}{\rm
Let $\cc T$ be a triangulated category. Then $\cc T$ is additive
and comes with an automorphism $\Sigma$. A square is defined to be
special iff its fold is a distinguished triangle in $\cc T$. When
we think of a triangulated category $\cc T$ as being the category
with squares defined above, then we shall denote it as $\cc T^d$.

Let $\cc A$ be an abelian category. Let $Gr^b\cc A$ be the
category of bounded, graded objects in $\cc A$. We recall the
reader that a graded object of $\cc A$ is a sequence of objects
$\{A_i\in\cc A\}_{i\in\bb Z}$. The sequence $\{A_i\}$ is bounded
if $A_i=0$ except for finitely many $i\in\bb Z$.

We define the functor $\Sigma:Gr^b\cc A\to Gr^b\cc A$ to be the
shift, that is $\Sigma\{A_i\}=\{B_i\}$ with $B_i=A_{i+1}$. A
square in $Gr^b\cc A$ is defined to be special if the fold
   $$A\xrightarrow{(\alpha,-\beta)^t}B\ps C\xrightarrow{(\gamma,\delta)}D\lra{\mu}\Sigma A$$
gives a long exact sequence in $\cc A$
   $$\cdots\to D_{i-1}\to A_i\to B_i\ps C_i\to D_i\to A_{i+1}\to\cdots$$

Let $H:D^b(\cc A)\to Gr^b\cc A$ be the homology functor taking a
complex $A\in D^b(\cc A)$ to $\{H_i(A)\}$. Then it induces a
functor between categories with squares
   $$H:D^b(\cc A)^d\to Gr^b\cc A.$$

}\end{exs}

\begin{defs}{\rm
(1) Let $\cc T$ be a category with squares and $m,n\geq 0$. A
functor $X:\Delta^m\times\Delta^n\to\cc T$ is called an {\it
augmented diagram\/} if for any $0\le i\le i'\le m$ and $0\le j\le
j'\le n$ we are given a special square
   $$\xymatrix{X_{ij'}\ar[r]&X_{i'j'}\ar@/^55pt/[dl]^{\delta^{i',j'}_{i,j}}\\
               X_{ij}\ar[u]\ar[r]&X_{i'j}\ar[u]}$$
such that $\delta^{i',j'}_{i,j}$ is the composite
   $$X_{i'j'}\to X_{mn}\xrightarrow{\delta_{0,0}^{m,n}}\Sigma X_{00}\to\Sigma X_{ij}.$$
By a {\it morphism\/} between augmented diagrams $\phi:X\to Y$ is
meant a natural transformation of functors such that the square
   $$\xymatrix{X_{i'j'}\ar[d]_{\phi_{i'j'}}\ar[r]^{\delta^{i',j'}_{i,j}}&\Sigma X_{ij}\ar[d]^{\Sigma\phi_{ij}}\\
               Y_{i'j'}\ar[r]^{\delta^{i',j'}_{i,j}}&\Sigma Y_{ij}}$$
is commutative for any $0\le i\le i'\le m$ and $0\le j\le j'\le
n$.

The category of augmented diagrams will be denoted by $Q_{m,n}\cc
T$. There results a bisimplicial category $Q\cc T=\{Q_{m,n}\cc
T\}_{m,n\geq 0}$ (the face/degeneracy operators are defined by
deleting/inserting a row or column).

(2) For a category with squares $\cc T$, its $K$-theory $K(\cc T)$
is defined to be the space $\Omega|\Ob(Q\cc T)|$.

}\end{defs}

Let $H:D^b(\cc A)^d\to Gr^b\cc A$ be the functor of categories
with squares constructed above. We have the map of bisimplicial
categories $\chi:QD^b(\cc A)^d\to QGr^b\cc A$ induced by $H$, and
hence the map $K(\chi):K(D^b(\cc A))\to K(Gr^b\cc A)$.

Let $\cc E$ be an exact category and $m,n\geq 0$. Denote by
$Q_{m,n}\cc E$ the following category. Its objects are the
functors $X:\Delta^m\times\Delta^n\to\cc E$ such that for any
$0\le i\le i'\le m$ and $0\le j\le j'\le n$ we are given a
bicartesian square
   $$\xymatrix{X_{ij'}\ar@{>->}[r]&X_{i'j'}\\
               X_{ij}\ar@{->>}[u]\ar@{>->}[r]&X_{i'j}\ar@{->>}[u]}$$
in which the vertical arrows are epimorphisms and the horizontal
arrows are monomorphisms. The morphisms are defined by natural
transformations. The resulting bisimplicial category denote by
$Q\cc E$. It is well-known that a simplicial model for a delooping
of the space $K(\cc E)$ is given by the realization of the
bisimplisial set $\Ob Q\cc E$.

Let $\cc A$ be an abelian category and let $i:\cc A\to D^b(\cc A)$
denote the natural functor sending an object $A\in\cc A$ to the
complex concentrated in the zeroth degree. Then it induces a
functor (see also some discussion below) of bisimplicial
categories $\iota:Q\cc A\to Q D^b(\cc A)^d$. Note that the
differentials $\delta^{i',j'}_{i,j}$ in $Q D^b(\cc A)^d$ are
canonically unique for every diagram coming from $Q\cc A$
(see~\cite{N1}).

\begin{thm}[Neeman~\cite{N}]\label{neem}
Let $\cc A$ be a small abelian category. Then the composite
   $$\Ob Q\cc A\lra{\iota}\Ob Q D^b(\cc A)^d\lra{\chi}\Ob QGr^b(\cc A)$$
is a homotopy equivalence.
\end{thm}

As usual, given a category $\cc C$ denote by $i\cc C$ the maximal
groupoid in $\cc C$ and by $i.\cc C$ the nerve in the
$i$-direction.

\begin{cor}\label{ta}
Let $\cc A$ be a small abelian category. Then the composite of
maps of trisimplicial objects
   $$i.Q\cc A\lra{\iota}i.Q D^b(\cc A)^d\lra{\chi}i.QGr^b(\cc A)$$
is a homotopy equivalence.
\end{cor}

\begin{proof}
Given $k\geq 0$ the category $i_k\cc A$ of strings of isomorphisms
$A_0\lra{\sim}\cdots\lra{\sim}A_k$ is abelian and the composite
   $$i_kQ\cc A=Q[i_k\cc A]\lra{\iota}i_kQ D^b(\cc A)^d\lra{\chi}i_kQGr^b(\cc A)=QGr^b[i_k\cc A]$$
is a homotopy equivalence of bisimplicial objects by
Theorem~\ref{neem}. It follows from Lemma~\ref{segal} that is so
the map of the corollary.
\end{proof}

\renewcommand{\proofname}{Proof of Theorem~\ref{neeman}}

\begin{proof}
(1) First prove the statement for an abelian category $\cc A$. For
Quillen's $K$-theory $K(\cc A)$ we use the following simplicial
model. It is the loop space of the realization of $i.Q\cc A$
(see~\cite{W}). In turn, the model for $K(\bb D^b(\cc A))$ is
given by the bisimplicial maximal groupoid $iQ\bb D^b(\cc A)$ (see
section~\ref{s}).

By Corollary~\ref{ta} it suffices to show that the map $i.Q\cc
A\lra{\iota}i.Q D^b(\cc A)^d$ factors through $i.Q\bb D^b(\cc A)$.
Recall that $Q_{m,n}\bb D^b(\cc A)$, $m,n\geq 0$, consists of the
objects $X\in\bb D^b(\cc A)_{\Delta^m\times\Delta^n}$ such that
for any $0\le i\le i'\le m$ and $0\le j\le j'\le n$ the square
   $$\xymatrix{X_{ij'}\ar[r]&X_{i'j'}\\
               X_{ij}\ar[u]\ar[r]&X_{i'j}\ar[u]}$$
is bicartesian in $\bb D^b(\cc A)_\Box$ (= cocartesian in
triangulated d\'erivateurs~\cite{M}). It follows that
   $$\cone(X_{ij}\to X_{ij'}\ps X_{i'j})\to\cone(0\to X_{i'j'})\iso X_{i'j'}$$
is a quasi-isomorphism in $C^b(\cc A)$, hence an isomorphism in
$D^b(\cc A)$ (we use here properties of triangulated d\'erivateurs
and the triangulated structure information which $\bb D^b(\cc A)$
encodes~\cite{C1,F,M1}). Now compose the inverse of this
isomorphism with the natural projection
   $$\cone(X_{ij}\to X_{ij'}\ps X_{i'j})\to\cone(X_{ij}\to 0)\iso\Sigma X_{ij}$$
and we have a map $\delta^{i',j'}_{i,j}:X_{i'j'}\to\Sigma X_{ij}$.
This produces a special square in $D^b(\cc A)^d$.

The construction is clearly natural. Let $f:X\to Y$ with $X,Y\in
Q_{m,n}\bb D^b(\cc A)$ be an isomorphism. It is represented by a
diagram
   $$X\longleftarrow Z\to Y$$
with $Z\in Q_{m,n}\bb D^b(\cc A)$ and arrows quasi-isomorphisms.
We have the following commutative diagram in $C^b(\cc A)$ for any
$0\le i\le i'\le m$ and $0\le j\le j'\le n$:
   $$\xymatrix{
     X_{ij}\ar[r]^(.35){a} & X_{ij'}\ps X_{i'j}\ar[r] & X_{i'j'}&\ar[l]\cone(a)\ar[r] & \Sigma X_{ij}\\
     Z_{ij}\ar[r]^(.35){b}\ar[u]\ar[d]&Z_{ij'}\ps Z_{i'j}\ar[r]\ar[u]\ar[d]
     &Z_{i'j'}\ar[u]\ar[d]&\ar[l]\cone(b)\ar[u]\ar[d]\ar[r] & \Sigma Z_{ij}\ar[u]\ar[d]\\
     Y_{ij}\ar[r]^(.35){c} & Y_{ij'}\ps Y_{i'j}\ar[r] & Y_{i'j'}&\ar[l]\cone(c)\ar[r] & \Sigma Y_{ij}
     }$$
This yields an isomorphism of triangles in $D^b(\cc A)$
   $$\xymatrix{
     X_{ij}\ar[r]\ar[d]&X_{ij'}\ps X_{i'j}\ar[r]\ar[d]&X_{i'j'}\ar[d]\ar[r]^{\delta^{i',j'}_{i,j}}&\Sigma X_{ij}\ar[d]\\
     Y_{ij}\ar[r]&Y_{ij'}\ps Y_{i'j}\ar[r]&Y_{i'j'}\ar[r]^{\delta^{i',j'}_{i,j}}&\Sigma Y_{ij}}$$
and hence an isomorphism of special squares in $D^b(\cc A)^d$.

Now, let $X\in Q_{m,n}\bb D^b(\cc A)$ and $0\le i'\le m$ and $0\le
j'\le n$. There is a commutative diagram in $C^b(\cc A)$
   $$\xymatrix{
     X_{00}\ar[r]\ar[d]&X_{i'0}\ps X_{0j'}\ar[r]\ar[d]&X_{i'j'}\ar[d]\\
     X_{00}\ar[r]&X_{m0}\ps X_{0n}\ar[r]&X_{mn}}$$
and hence, in $D^b(\cc A)$ we deduce a commutative square
   $$\xymatrix{
     X_{i'j'}\ar[d]\ar[r]^{\delta^{i',j'}_{0,0}}&\Sigma X_{00}\ar[d]\\
     X_{mn}\ar[r]\ar[r]^{\delta^{m,n}_{0,0}}&\Sigma X_{00}}$$
Given any $0\le i\le i'$ and $0\le j\le j'$ there is a commutative
diagram in $C^b(\cc A)$
   $$\xymatrix{
     X_{00}\ar[r]\ar[d]&X_{i'0}\ps X_{0j'}\ar[r]\ar[d]&X_{i'j'}\ar[d]\\
     X_{ij}\ar[r]&X_{i'j}\ps X_{ij'}\ar[r]&X_{i'j'}}$$
and hence, in $D^b(\cc A)$ we deduce a commutative square
   $$\xymatrix{
     X_{i'j'}\ar[d]\ar[r]^{\delta^{i',j'}_{0,0}}&\Sigma X_{00}\ar[d]\\
     X_{i'j'}\ar[r]\ar[r]^{\delta^{i',j'}_{i,j}}&\Sigma X_{ij}}$$
and hence the ``natural'' map $\delta^{i',j'}_{i,j}:X_{i'j'}\to
\Sigma X_{ij}$ is obtained from $\delta^{m,n}_{0,0}:X_{mn}\to
\Sigma X_{00}$ just as the composite
   $$X_{i'j'}\to X_{mn}\xrightarrow{\delta^{m,n}_{0,0}}\Sigma X_{00}\to\Sigma X_{ij}.$$

It follows that the functors $dia:\bb D^b(\cc
A)_{\Delta^m\times\Delta^n}\to\Hom(\Delta^m\times\Delta^n,D^b(\cc
A))$, $m,n\geq 0$, induce a map of bisimplicial groupoids
   $$\theta:iQ\bb D^b(\cc A)\to iQ D^b(\cc A)^d.$$
Obviously, the map $i.Q\cc A\lra{\iota}i.Q D^b(\cc A)^d$ factors
as
   $$i.Q\cc A\lra{\rho}i.Q\bb D^b(\cc A)\lra{\theta}i.Q D^b(\cc A)^d.$$
This implies the claim.

(2) Suppose now that an exact category $\cc E\subseteq\cc A$
satisfies the assumptions of the theorem. Consider the commutative
diagram
   $$\xymatrix{i.Q\cc E\ar[d]\ar[r]^(.4)\rho & i.Q\bb D^b(\cc E)\ar[d]\\
              i.Q\cc A\ar[r]^(.4)\rho & i.Q\bb D^b(\cc A)\ar[r]^{\chi\theta}&i.QGr^b\cc A}$$
in which the vertical arrows are induced by the inclusion $\cc
E\to\cc A$. The left vertical arrow is a homotopy equivalence by
the Resolution Theorem~\cite{Q}. The fact that the map
$\chi\theta\rho$ is a homotopy equivalence by~(1) obviously
finishes the proof.
\end{proof}

\renewcommand{\proofname}{Proof}

Basing on Vaknin's computations~\cite{V} Neeman
shows~\cite[p.~39]{N1} that there is an exact category $\cc E$
such that the homomorphism $K_1(\iota):K_1(\cc E)\to K_1(D^b(\cc
E))$ is not a monomorphism (while it is a split monomorphism for
abelian categories~\cite{N,N1}). The simplest example is where
$\cc E$ is the category of free modules of finite rank over the
ring of dual numbers $k[\epsilon]/\epsilon^2$. Such exact
categories could give us counter-examples to the first
Maltsiniotis conjecture if we showed in a similar way that the map
$K_1(\rho):K_1(\cc E)\to K_1(\bb D^b(\cc E))$ is not a
monomorphism.

\end{document}